\documentclass[11pt]{amsart}
\usepackage{amssymb}
\usepackage{amsmath, amsthm, amsopn, amsfonts}
\usepackage[dvips]{graphicx}
\usepackage{graphpap}
\newtheorem{theoreme}{Theorem}

\newtheorem{proposition}{Proposition}
\newtheorem{lemme}[proposition]{Lemma}

\newtheorem{remarque}[proposition]{Remark}

\numberwithin{equation}{section}
\numberwithin{proposition}{section}

\def\Im{\textrm{Im}}
 
\def\11{{\rm 1~\hspace{-1.4ex}l} }
\def\R{\mathbb R}
\def\C{\mathbb C}
\def\Z{\mathbb Z}
\def\N{\mathbb N}

\def\T{\mathbb T}

\begin{document}
\title[On Strichartz estimates]
{Strichartz estimates for long range perturbations }
\author{Jean-Marc Bouclet}
\address{D\'epartement de Math\'ematiques, Universit\'e Lille I, 59 655
Villeneuve d'Ascq Cedex, France}
\email{jean-marc.bouclet@math.univ-lille1.fr}
\author{Nikolay Tzvetkov}
\address{D\'epartement de Math\'ematiques, Universit\'e Lille I, 59 655
Villeneuve d'Ascq Cedex, France}
\email{nikolay.tzvetkov@math.univ-lille1.fr}
\begin{abstract} 
We study local in time Strichartz estimates for the Schr\"odinger equation
associated to long range perturbations of the flat Laplacian on the euclidean
space. We prove that in such a geometric situation, outside a large ball
centered at the origin, the solutions of the Schr\"odinger equation enjoy the
same Strichartz estimates as in the non perturbed situation. The proof is based
on the Isozaki-Kitada parametrix construction. If in addition the metric is
non trapping, we prove that the Strichartz estimates hold in the whole space.
\end{abstract}
\maketitle
\tableofcontents

\section{Introduction.}
Let $(M,g)$ be a $d$-dimensional Riemannian manifold. 
Denote by $\Delta_{g}$ the Laplace-Beltrami operator associated to the metric
$g$. Consider the time dependent Schr\"odinger equation on $(M,g)$
\begin{equation}\label{1}
iu_{t}+\Delta_{g}u=0
\end{equation}
subject to initial data
\begin{equation}\label{2}
u|_{t=0}=u_{0}\, .
\end{equation}
It is well-known (see e.g. \cite{GV}) that when $(M,g)$ is the
flat Euclidean space 
(i.e. $\R^d$ with the metric $g_{ij}=\delta_{ij}$ the unit $d\times d$ matrix) the
solutions of (\ref{1})-(\ref{2}) enjoy the (local in time) Strichartz
estimates
\begin{equation}\label{3} 
\|u\|_{L^{p}([0,1];L^{q}(\R^d))}\leq C\|u_0\|_{L^2(\R^d)}\,\, ,
\end{equation}
where
\begin{equation}\label{adm}
\frac{2}{p}+\frac{d}{q}=\frac{d}{2},\quad p\geq 2,\quad (p,q)\neq (2,\infty)\, .
\end{equation}
Moreover (\ref{3}) is global in time which means that one can replace $[0,1]$
in the left hand-side of (\ref{3}) by $\R$. In \cite{Bo,BGT1} one studies the
possible extensions of (\ref{3}) to the situation where $M$ is compact. An important
new phenomenon that one has to take into account, when $M$ is compact, is the 
unavoidable derivative loss in (\ref{3}) for some values of $(p,q)$. By
``loss'' we mean that $\|u_0\|_{L^2}$ in the right hand-side of (\ref{3})
{\it should} be replaced by $\|u_0\|_{H^s}$ for some positive $s$. Here are
two significant examples. If $M$ is the standard sphere $S^d$, $d\geq 3$,
then it is proved in \cite{BGT1} that the solutions of (\ref{1})-(\ref{2})
satisfy
\begin{equation}\label{4}
\|u\|_{L^{2}\big([0,1];L^{\frac{2d}{d-2}}(S^d)\big)}\leq C\|u_0\|_{H^{\frac{1}{2}}(S^d)}
\end{equation}
(notice that the couple $(2,\frac{2d}{d-2})$ satisfies (\ref{adm})). Moreover,
the $H^{1/2}(S^d)$ norm in the right hand-side of (\ref{4}) is sharp in the
sense that for every $\varepsilon>0$, the estimate
\begin{equation*}
\|u\|_{L^{2}\big([0,1];L^{\frac{2d}{d-2}}(S^d)\big)}\leq C\|u_0\|_{H^{\frac{1}{2}-\varepsilon}(S^d)}
\end{equation*}
is false. 

A second example where one should encounter losses in (\ref{3}) is the flat torus
$\T^d=\R^d| (2\pi\Z)^{d}$. 
More precisely the estimate
\begin{equation}\label{5}
\|u\|_{L^{\frac{2(d+2)}{d}}\big([0,1];L^{\frac{2(d+2)}{d}}(\T^d)\big)}\leq C\|u_0\|_{L^{2}(\T^d)}
\end{equation}
is false (notice that again the couple $(\frac{2(d+2)}{d},\frac{2(d+2)}{d})$ satisfies (\ref{adm})).
We refer to \cite{Bo} for a counterexample disproving (\ref{5}) in the case
$d=1$. The extension to higher dimensions is straightforward. One may however
expect (\ref{5}) to be replaced by
\begin{equation}\label{6}
\|u\|_{L^{\frac{2(d+2)}{d}}\big([0,1];L^{\frac{2(d+2)}{d}}(\T^d)\big)}\leq
C_{s}\|u_0\|_{H^{s}(\T^d)},\quad s>0.
\end{equation}
Estimate (\ref{6}) is known for $d=1,2$ (see \cite{Bo}) (in this case
$\frac{2(d+2)}{d}$ is an even integer). For $d\geq 3$, the study of (\ref{6})
leads to an interesting open problem.
\\

When $M$ is not compact, extensions of (\ref{3}) 
were recently studied by several authors (see \cite{ST,BGT1,RZ,HTW}). 
In the works \cite{ST,RZ,HTW} the authors consider
non compact manifolds with metrics which are a ``small'' perturbation at infinity of a fixed
``nice'' metric, satisfying a non trapping assumption on the geodesic
flow. It turns out that in such a geometric situation, one can prove exactly the
same estimates as for the Euclidean space. 

In \cite{BGT1}, one considers
$(M,g)$ to be $\R^d$ with a perturbation of the flat metric {\it without} the non trapping assumption.
In this context
one can get the Strichartz estimates with losses, just as in the case of a
compact manifold. It is however a priori not clear whether losses of derivatives
in the Strichartz estimates may come from the geometry at infinity. The main goal of this paper is to show
that one {\it can not} have losses in the Strichartz inequalities coming from the
geometry at infinity in the case of long range perturbations of the euclidean
metric on $\R^d$.
\begin{theoreme}\label{thm1}
Consider $\R^d$ equipped with a smooth Riemannian metric 
$$
g (x)  =(g_{ij}(x))_{i,j=1}^{d},\quad x\in\R^d
$$ 
satisfying for some
$\rho>0$,
\begin{equation}\label{H1}
|\partial^{\alpha}(g^{ij}(x)-\delta_{ij})|\leq C_{\alpha}\langle
 x\rangle^{-\rho-|\alpha|},\quad i,j\in\{1,\dots d\}
\end{equation}
($\delta_{ij}$ being the Kronecker symbol) and
\begin{equation}\label{H1bis}
c\,{\rm Id}\leq g(x) \leq C\, {\rm Id}\, .
\end{equation}
Then there exists $R>0$ such that for every $T>0$, every $(p,q)$ satisfying
(\ref{adm}) there exists $C>0$ 
such that for every $f\in L^2(\R^d)$,
\begin{equation}\label{est1}
\|e^{it\Delta_{g}}f\|_{L^{p}([-T,T];L^{q}(|x|\geq R))}\leq
C\|f\|_{L^2(\R^d)}\, .
\end{equation}
Moreover, for every $f\in H^{\frac{1}{p}}(\R^d)$,
\begin{equation}\label{est2}
\|e^{it\Delta_{g}}f\|_{L^{p}([-T,T];L^{q}(|x|\leq R))}\leq
C\|f\|_{H^{\frac{1}{p}}(\R^d)}\, .
\end{equation}
\end{theoreme}
\begin{remarque}
The result of Theorem~\ref{thm1} is stated only for metric perturbations of
the flat Laplacian. However, an examination of our proof shows that the
statement still holds if we add long range lower order terms. The same remark
is valid for  Theorem~\ref{thm2} below.
\end{remarque}
\begin{remarque}
Let us emphasize that estimates in the spirit of
(\ref{est1}) are known to hold in the context of resolvent estimates for
long range perturbations of the Laplacian (cf. \cite{burqreso,vodev}).
By this we mean the following: the fact that we have no derivative loss in (\ref{est1})
as in the free case is somehow similar to the fact that the high energy
resolvent estimates with weights supported near infinity 
are the same as for the free resolvent.
\end{remarque}

If we suppose that the metric $g$ is non trapping then one can improve
(\ref{est2}) and get the full family of Strichartz estimates. Recall that $g$
is non trapping if every geodesic (globally defined thanks to (\ref{H1bis}))
leaves every compact set in finite time. Let us now state our second result.
\begin{theoreme}\label{thm2}
Under the assumptions of Theorem~\ref{thm1}, if in addition we suppose that $g$
is non trapping, then 
\begin{equation}\label{est3}
\|e^{it\Delta_{g}}f\|_{L^{p}([-T,T];L^{q}(\R^d))}\leq C\|f\|_{L^2(\R^d)}\,.
\end{equation}
\end{theoreme}
Note that under the short range condition $\rho>1$,
estimate (\ref{est3}) is proved in
\cite{RZ} by using FBI transform techniques.

Let us now explain the main points in the proof of the above results. The
proof of (\ref{est1}) is based on the Isozaki-Kitada \cite{IK} parametrix
construction. Recall that this construction was introduced to build 
modified scattering operators for long range perturbations of the free
Schr\"odinger group. 
Let us point out that,
since here we are only dealing with finite time estimates,
we are not using the Isozaki-Kitada method in its full strength. 
In particular, we do not need a non trapping assumption on the metric to get
(\ref{est1}). If we were interested in proving (\ref{est1}) with a constant $C$
uniform with respect to $T$, then a non trapping assumption and the full force
of the  Isozaki-Kitada method would be needed. We will not address this
interesting issue here. 
See \cite{burq,RT} for the proof of the global in time estimates in the case
of compactly supported perturbations.
The proof of (\ref{est2}) is essentially contained in
\cite{BGT1}. 
The proof of Theorem~\ref{thm2} is based on ideas introduced in \cite{ST,BGT1}.
In fact, it is fair to say that, as far as the spatial regularity is concerned,
the estimates established in \cite{BGT1} are gaining $1/2$ derivative with
respect to the Sobolev embedding. 
We prove Theorem~\ref{thm2} by showing
that the missing $1/2$ derivative can be recovered thanks to the local
smoothing effect (when it is available). 
Let us notice that this effect is a consequence of standard resolvent
estimates for non trapping perturbations of the Laplacian.
It would be interesting to know whether
intermediate situations may exist and if so to quantify them in terms of the metric. 
It is worth mentioning the work \cite{burq}, where (\ref{est3}) with
$C_{\varepsilon}\|f\|_{H^{\varepsilon}}$, $\varepsilon>0$ 
instead of $C\|f\|_{L^2}$ is studied, i.e. an unnecessary $\varepsilon$ derivative
loss is accepted. In this context, let us recall that the analysis in \cite{BGT1,BGT2} has shown
that, if one is interested in non linear applications, the losses in term of
Sobolev regularity in the right hand-side of (\ref{est3}) are more dramatic
then the losses in terms of the range of
possible values of $(p,q)$ in the left hand-side of (\ref{est3}). 

The rest of this paper is organized as follows. In the next section, we
introduce the functional calculus for $\Delta_g$ in terms of pseudo
differential calculus. In section 3, we recall the main points of the
Isozaki-Kitada parametrix. The analysis of section 3 is then used in section 4
for the proof of (\ref{est1}). Section 5 deals with estimates on time
intervals depending on the frequency localization. In section 6, we use the
non trapping condition to get the full family of Strichartz estimates in a fixed
compact set. Section 7 is devoted to the rather standard non linear applications of the
analysis of the previous sections.

{\bf Notation.} In this paper several numerical constants will be denoted by
the same $C$. For $T>0$, $p\in[1,+\infty]$, and $B$ a Banach space, we denote by $L^p_TB$ the
Banach space of $L^p$  functions on $[-T,T]$ with values in $B$ equipped with the natural norm.
We denote by $P\geq 0$ the self adjoint realization of $-\Delta_{g}$ on $L^2(\R^d)$.
\section{Functional calculus}
In this section $g$ is a metric on $\R^d$ such that there exist $C\geq 1\geq
c>0$ such that
\begin{equation}\label{h1}
c|\xi|^{2}\leq g(x)(\xi,\xi)\leq C|\xi|^{2},\quad \forall x\in \R^d,\,\,
\forall \xi\in\R^d\,,
\end{equation}
and for every $\alpha\in\N^d$,
\begin{equation}\label{h2}
|\partial_{x}^{\alpha}\, g(x)|\leq C_{\alpha}\,\, .
\end{equation}
Notice that we do not assume the long range condition (\ref{H1}).

Further, we denote by
$$
p_0(x,\xi):=g(x)(\xi,\xi)\equiv\sum_{i,j=1}^{d}g^{ij}(x)\xi_i\xi_j
$$
the principal symbol of $-\Delta_g$.
Here we adopt the standard notation for $(g^{ij}(x))\equiv(g_{ij}(x))^{-1}$.

The goal of this section is 
to approximate  $\varphi(-h^2 \Delta_{g})$,
$h\in]0,1]$,
for a given bump function
$\varphi\,:\, \R\rightarrow\R$,  by a suitable semi-classical pseudo differential operator.
Similar considerations are performed in \cite{BGT1}, where $\Delta_g$ is the
Laplace-Beltrami operator on a {\it compact} Riemannian manifold. Here we
follow a similar scheme. The new feature in our analysis is the $L^p$ bound
for the remainder of the pseudo differential expansion of $\varphi(-h^2 \Delta_{g})$.
In \cite{BGT1} this is done
by only invoking $L^2$ considerations and the fact that, on a compact manifold
$M$, $L^{\infty}(M)$ is continuously embedded
in $L^{2}(M)$. Here, we can not use this fact. We overcome the difficulty
by using $L^p$ bounds for powers of the resolvent of $P$.

Recall that $P\geq 0$ is the self adjoint realization of $-\Delta_{g}$ on
$L^2(\R^d)$. We first collect several classical properties of $P$. For every
$s\in\R$, there exists a constant $C_s$ such that for every 
$u\in {\mathcal S}(\R^d)$,
\begin{equation}\label{lesno}
C_{s}^{-1}\|(P+1)^{s/2}u\|_{L^2(\R^d)}\leq \|u\|_{H^s(\R^d)}\leq
C_{s}\|(P+1)^{s/2}u\|_{L^2(\R^d)}\,.
\end{equation}
Next, we recall that the Schwartz class is stable under the action of the
resolvent of $P$. More precisely, for every $ z\in\C\backslash[0,\infty[$
the map $(P-z)^{-1}$ is acting continuously on ${\mathcal S}(\R^d)$.
As a consequence, by the standard duality argument it acts continuously on 
${\mathcal S}'(\R^d)$ too. In particular, for every $w\in L^p(\R^d)$, $1\leq
p\leq +\infty$, $(P-z)^{-1}w$ is well defined.

The elliptic nature of $P$ also implies that, for every $s\in\R$, there exists
$C_s$ such that for every $u\in {\mathcal S}(\R^d)$,
$$
\|(P+1)^{-1}u\|_{H^{s+2}(\R^d)}\leq C_s\|u\|_{H^s(\R^d)}.
$$
As a consequence for every $z\in\C\backslash[0,\infty[$,
\begin{equation}\label{im}
\|(P-z)^{-1}u\|_{H^{s+2}(\R^d)}\leq \frac{\widetilde{C}_s\langle z\rangle}{|\Im\, z|}\,
\|u\|_{H^s(\R^d)}.
\end{equation}
Indeed, we can write
\begin{eqnarray*}
\|(P-z)^{-1}u\|_{H^{s+2}} & \leq & C\|(P+1)^{\frac{s+2}{2}}(P-z)^{-1}u\|_{L^2}
\\
& \leq & C\|(P-z)^{-1}(P+1) (P+1)^{\frac{s+2}{2}}(P+1)^{-1}u\|_{L^2}
\\
&\leq &
C\|(P-z)^{-1}(P+1)\|_{L^2\rightarrow L^2 }\|(P+1)^{\frac{s}{2}}u\|_{L^2}
\\
& \leq &
C\,\sup_{\lambda\in \R^{+}}\frac{|\lambda+1|}{|\lambda-z|}\|u\|_{H^s}
\\
&\leq &
\frac{C\langle z\rangle}{|\Im\, z|}\,
\|u\|_{H^s}
\end{eqnarray*}
which proves (\ref{im}).

We next state a bound for sufficiently large powers of the resolvent of $P$.
\begin{proposition}\label{neg}
Let $n_0>d/2$ be an integer. For $z\in \C\backslash[0,\infty[$, we denote by
$K_{z}(x,y)$ the Schwartz kernel of the operator $(P-z)^{-n_0}$. Then for
every $\alpha\in\N^d$, there exist $C_{\alpha}>0$ and $n(\alpha)\in\N$ such
that for every $(x,y)\in\R^{2d}$, and every $z\in \C\backslash[0,\infty[$,
\begin{equation*}
\big|(x-y)^{\alpha}K_{z}(x,y)\big|\leq C_{\alpha}
\Big(\frac{\langle z\rangle}{|{\rm Im} \, z|}\Big)^{n(\alpha)}\, .
\end{equation*}
In particular there exist $N$ and $C>0$ such that
$$
|K_{z}(x,y)|\leq C\,\Big(\frac{\langle z\rangle}{|{\rm Im}\, z|}\Big)^{N}\,
\frac{1}{(1+|x-y|^2)^{d}}\,\, .
$$
\end{proposition}
\begin{proof}
Let first $|\alpha|=0$.
Observe that for $s>d/2$ any bounded linear map
$$
A\,:\, H^{-s}(\R^d)\longrightarrow H^s(\R^d)
$$
has a Schwartz kernel $K_{A}(x,y)$ which is a bounded continuous function on
$\R^{2d}$.
Moreover, there exists $C>0$ such that for every bounded map from $H^{-s}$ to
$H^s$, every $(x,y)\in\R^{2d}$,
\begin{equation}\label{Dirak}
|K_{A}(x,y)|\leq C\||A\||,
\end{equation}
where $\||A\||$ denotes the norm of $A$,
$$
\||A\||:=\sup_{\|u\|_{H^{-s}=1}}\|Au\|_{H^s}\,.
$$
Using (\ref{im}), we get the estimate
$$
\|(P-z)^{-n_0}u\|_{H^{n_0}}\leq C\Big(\frac{\langle z\rangle}{|\Im\, z|}\Big)^{n_0}\|u\|_{H^{-n_0}}\, .
$$
Therefore, in view of (\ref{Dirak}), the assertion of
Proposition~\ref{neg} holds for $|\alpha|=0$.

Let next $|\alpha|=1$.
Clearly, for $k\in\{1,\dots,d\}$, $(x_k-y_k)K_{z}(x,y)$ is the
kernel of $[x_k,(P-z)^{-n_0}]$. Using the identity
$$
\big[x_k,(P-z)^{-1}\big]=(P-z)^{-1}[P,x_k](P-z)^{-1},
$$
we arrive at the formula
\begin{eqnarray*}
\big[x_k,(P-z)^{-n_0}\big] & = &
\sum_{j=0}^{n_0-1}(P-z)^{-(j+1)}[P,x_k](P-z)^{-(n_0-j)}
\\
& = &
\sum_{j_1+j_2=n_0-1}
(P-z)^{-(j_1+1)}[P,x_k](P-z)^{-(j_2+1)}\, .
\end{eqnarray*}
On the other hand, by invoking (\ref{im}), we obtain that if $j_1+j_2=n_0-1$,
and if $s\in\R$ is such that $n_0>s>d/2$, then
the linear map
$$
(P-z)^{-(j_1+1)}[P,x_k](P-z)^{-(j_2+1)}
$$
is bounded from $H^{-s}(\R^d)$ to 
$$
H^{-s+2(j_2+1)-1+2(j_1+1)}(\R^d)=H^{-s+2n_0+1}(\R^d)\subset H^{s}(\R^d)
$$
with operator norm bounded by a polynomial of 
\begin{equation}\label{frac}
\frac{\langle z\rangle}{|\Im\, z|}
\end{equation}
which, using (\ref{Dirak}), yields the assertion of
Proposition~\ref{neg} for $|\alpha|=1$.

Let finally $|\alpha|$ be arbitrary.
For $k\in\{1,\dots,d\}$, $\alpha\in\N^d$ and $A$ a function of $P$, we introduce the notations
$$
\zeta_{k}A:=[x_k,A],\quad
\zeta^{\alpha}A=\zeta_{1}^{\alpha_1}\dots\zeta_{d}^{\alpha_d}A\, .
$$
Notice that $\zeta_{k}\zeta_{j}A=\zeta_{j}\zeta_{k}A$. Using an induction
argument, one can check that
$$
\zeta^{\alpha}(P-z)^{-n_0}
$$
is a linear combination of terms of type
$$
(P-z)^{-(1+j_1)}(\zeta^{\alpha^1}P)(P-z)^{-(1+j_2)}
(\zeta^{\alpha^2}P)\dots
(P-z)^{-(1+j_n)}(\zeta^{\alpha^n}P)(P-z)^{-(1+j_{n+1})}
$$
with
$$
\alpha^{i}\neq 0,\quad
1\leq n\leq |\alpha|,\quad
\alpha^1+\dots +\alpha^{n}=\alpha,\quad
j_1+j_2+\dots + j_{n+1}=n_0-1\, .
$$
Therefore for $n_0>s>d/2$, the map
$
\zeta^{\alpha}(P-z)^{-n_0}
$
is bounded map from $H^{-s}(\R^d)$ to $H^{2n_0+|\alpha|-s}(\R^d)$
with operator norm bounded by a polynomial of (\ref{frac}).
By invoking (\ref{Dirak}), we complete the proof of Proposition~\ref{neg}.
\end{proof}
We will use the result of Proposition~\ref{neg} to get $L^p$ bounds for  
sufficiently large powers of the resolvent of $P$.
For that purpose, we recall the well known Schur lemma.
\begin{proposition}\label{shur}
Let $K(x,y)$ be a continuous function on $\R^{2d}$ satisfying
$$
\int_{\R^d}|K(x,y)|dx\leq C,\quad  \int_{\R^d}|K(x,y)|dy\leq C
$$
for some positive constant $C$. Then for every $p\in[1,+\infty]$, the linear map with kernel $K$ is bounded
on $L^p(\R^d)$ with norm $\leq C$.
\end{proposition}
\begin{proof}
The statement is straightforward for  $p=\infty$ and $p=1$. 
The case of an arbitrary $p$ follows then by interpolation 
(one can also easily give a direct proof avoiding the interpolation). 
\end{proof}
A direct combination of  Proposition~\ref{neg} and Proposition~\ref{shur}
gives the following statement.
\begin{proposition}\label{ppp}
Let us fix an integer $n_0>d/2$. 
Then there exist $C>0$ and $n\in \N$ such that for every $p\in [1,\infty]$,
$u\in L^p(\R^d)$, every $z\in \C\backslash[0,\infty[$,
$$
\|(P-z)^{-n_0}u\|_{L^p(\R^d)}
\leq C
\Big(\frac{\langle z\rangle}{|{\rm Im}\, z|}\Big)^{n}\|u\|_{L^p(\R^d)}\,.
$$
\end{proposition}
We next give another consequence of the Schur lemma.
\begin{proposition}\label{sympa}
Let $m>d$ and let $a(x,\xi,h)$ be a continuous function on $\R^d\times\R^d\times ]0,1]$
smooth with respect to the second variable satisfying
\begin{equation}\label{sab1}
\big|\partial_{\xi}^{\beta}a(x,\xi,h)\big| 
\leq 
C_{\beta}\langle\xi\rangle^{-m},\quad \beta\in\N^d \, .
\end{equation}
Then for $\infty \geq r\geq q\geq 1$,
\begin{equation}\label{sab2}
\|a(x,hD,h)\|_{L^q(\R^d)\rightarrow L^r(\R^d)}\leq C_{qr}\,
h^{d(\frac{1}{r}-\frac{1}{q})},\quad \forall h\in ]0,1]
\end{equation}
with
$$
C_{qr}\lesssim \sup_{|\beta|\leq d+1}C_{\beta}\,.
$$
\end{proposition}
\begin{proof}
Write
\begin{equation}\label{sab3}
a(x,hD,h)u(x)=\int_{\R^d}K_{h}(x,y)u(y)dy
\end{equation}
with
$$
K_{h}(x,y)=(2\pi h)^{-d}\int_{\R^d}e^{i\frac{(x-y)\cdot\xi}{h}}a(x,\xi,h)d\xi\, .
$$
Notice that, thanks to (\ref{sab1}) and the assumption $m>d$, the last integral
is absolutely convergent. Let us denote by $\hat{a}$ the Fourier transform of
$a$ with respect to the second variable. Then we can write
$$
K_{h}(x,y)=(2\pi h)^{-d}\,\hat{a}\big(x,\frac{x-y}{h},h\big)\,.
$$
Thanks to (\ref{sab1}), for every $\beta\in\N^d$,
$$
\big|
z^{\beta}\hat{a}(x,z,h)
\big|\lesssim C_{\beta}\, .
$$
Hence
$$
\sup_{x\in\R^d}\sup_{h\in]0,1]}\big|
\hat{a}(x,z,h)
\big|
\lesssim \langle z\rangle^{-d-1}\sup_{|\beta|\leq d+1}C_{\beta}\, .
$$
Therefore
$$
\int_{\R^d}|K_h(x,y)|dx
=
(2\pi)^{-d}
\int_{\R^d}|\hat{a}(y+hz,z,h)|dz\lesssim 
\sup_{|\beta|\leq d+1}C_{\beta}
$$
and
$$
\int_{\R^d}|K_h(x,y)|dy
=
(2\pi)^{-d}
\int_{\R^d}|\hat{a}(x,z,h)|dz\lesssim 
\sup_{|\beta|\leq d+1}C_{\beta}\, .
$$
Applying Proposition~\ref{shur} completes the proof of (\ref{sab2}) for $q=r$.
Thanks to (\ref{sab1}),
$$
|K_h(x,y)|\lesssim C_{0}\,h^{-d}\, ,
$$
and coming back to (\ref{sab3}), this completes the proof  of (\ref{sab2}) for
$r=\infty$ and $q=1$. Let us finally fix arbitrary $r,q$ satisfying 
$\infty\geq r\geq q\geq 1$. Interpolating between the $L^1\rightarrow L^1$ and
$L^1\rightarrow L^{\infty}$ bounds, we get the $L^1\rightarrow L^{r/q}$ bound.
Then, interpolation between the $L^{\infty}\rightarrow L^{\infty}$ and the
$L^1\rightarrow L^{r/q}$ bounds, we get the $L^{q}\rightarrow L^r$ bound.
This completes the proof of Proposition~\ref{sympa}.
\end{proof}
The next statement describes functions of $P$ in terms of semi-classical
pseudo differential operators.
Recall that $p_0$ denotes the principal symbol of $P$.
\begin{proposition}\label{cf}
Let $g$ be a metric on $\R^d$ satisfying (\ref{h1}) and (\ref{h2}).
Then for every $\varphi\in C_{0}^{\infty}(\R)$, there exist symbols
$(a_k)_{k\geq 0}$ satisfying
\begin{equation}\label{decay}
\big|
\partial_{x}^{\alpha}\partial_{\xi}^{\beta}a_{k}(x,\xi)
\big| 
\leq 
C_{k\alpha\beta},\quad \forall\, (x,\xi)\in {\rm supp}\,\varphi\circ p_{0},
\end{equation}
\begin{equation}\label{zero}
a_{k}(x,\xi)=0,\quad\forall\, (x,\xi)\notin {\rm supp}\,\varphi\circ p_{0},
\end{equation}
and there exists $n_1 \in\N$ such that for every $N\geq 1$, every $1\leq r\leq
q\leq\infty$, there exists $C_{Nqr}$ such that for every $h\in]0,1]$,
$$
\Big\|\varphi(h^2 P)-\sum_{k=0}^{N}h^ka_{k}(x,hD)\Big\|_{L^q(\R^d)\rightarrow L^r(\R^d)}\leq C_{Nqr}\,
h^{N-n_1+d(\frac{1}{r}-\frac{1}{q})}\, .
$$
Moreover for every $s\geq 0$ there exists $N_s$ such that for $N\geq N_s$,
$$
\Big\|\varphi(h^2 P)-\sum_{k=0}^{N}h^ka_{k}(x,hD)\Big\|_{H^{-s}(\R^d)\rightarrow H^s(\R^d)}\leq C\, h^{N-N_s}\, .
$$
\end{proposition}
\begin{remarque}
If we suppose that the metric $g$ satisfies (\ref{H1}) then we can replace the bound
(\ref{decay}) in Proposition~\ref{cf} by
\begin{equation}\label{decay-bis}
\big|
\partial_{x}^{\alpha}\partial_{\xi}^{\beta}a_{k}(x,\xi)
\big| 
\leq 
C_{k\alpha\beta}\langle x\rangle^{-k-|\alpha|},
\end{equation}
for every $(x,\xi)\in\R^{2d}$ such that $p_{0}(x,\xi)\in {\rm supp}\,\varphi$.
\end{remarque}
\begin{proof}[Proof of Proposition~\ref{cf}]
We first describe the classical construction (see
\cite{Seely,Robert,BGT1}) of a parametrix
for $(h^2P-z)^{-1}$, 
$z\in \C\backslash[0,\infty[$, $h\in]0,1]$. There exist
symbols 
$$
q_{0}(x,\xi,z),q_{1}(x,\xi,z), \dots,q_{N}(x,\xi,z)
$$
satisfying 
$$
\big|
\partial_{x}^{\alpha}\partial_{\xi}^{\beta}q_{k}(x,\xi,z)
\big| 
\leq 
C_{k\alpha\beta}\langle\xi\rangle^{-2-k-|\beta|}
\Big(\frac{\langle z\rangle}{|\Im\, z|}\Big)^{n(k,\alpha,\beta)}
$$
such that for every $h\in]0,1]$, and for every $N\geq 1$,
$$
(h^2P-z)\,\,\sum_{k=0}^{N}q_{k}(x,hD,z)={\rm Id }\, +h^{N+1}r_{N+1}(x,hD,z,h)
$$
with
$$
\big|
\partial_{x}^{\alpha}\partial_{\xi}^{\beta}r_{N+1}(x,\xi,z,h)
\big| 
\leq 
C_{N\alpha\beta}\langle\xi\rangle^{-N-1-|\beta|}
\Big(\frac{\langle z\rangle}{|\Im\, z|}\Big)^{n(N,\alpha,\beta)}
$$
uniformly in $h\in]0,1]$. 
Moreover the symbols are analytic with respect to $z\in
\C\backslash[0,\infty[$ and we can write
$$
q_{0}(x,\xi,z)=\frac{1}{p_0(x,\xi)-z}
$$
where $p_0(x,\xi)=g(x)(\xi,\xi)$ is the principal symbol of $-\Delta_{g}$. In
addition, for $k\geq 1$, $q_k$ takes the form
$$
q_{k}(x,\xi,z)=\sum_{j=1}^{2k-1}\,\frac{d_{j,k}(x,\xi)}{(p_0(x,\xi)-z)^{1+j}}
$$
with $d_{j,k}\in S^{2j-k}(\R^{2d})$ (they are polynomials in $\xi$ with coefficients
which are linear combinations of products
of derivatives of the coefficients of the inverse of the metric).

Therefore for every $h\in]0,1]$ every $z\in \C\backslash[0,\infty[$, every
$N\geq 1$,
$$
(h^2P-z)^{-1}=\sum_{k=0}^{N}h^kq_{k}(x,hD,z)-(h^2P-z)^{-1}h^{N+1}r_{N+1}(x,hD,z,h)\,.
$$
Then, for every $\varphi\in
C_{0}^{\infty}(\R)$, we can use the Helffer-Sj\"ostrand formula (see \cite{HS,DS}), 
$$
\varphi(h^2 P)=-\frac{1}{\pi}
\int_{\C}\bar{\partial}\widetilde{\varphi}(z)
(h^2 P-z)^{-1}dL(z),
$$
where $dL(z)$ denotes the Lebesgue measure on $\C$ and
$\widetilde{\varphi}(z)\in C_{0}^{\infty}(\C)$ 
is an almost analytic extension of
$\varphi$ which satisfies
\begin{equation}\label{anal}
\forall\, \Lambda>0,\quad |\bar{\partial}\widetilde{\varphi}(z)|\leq
C_{\Lambda}|\Im z|^{\Lambda}\, .
\end{equation}
This implies that $\varphi(h^2 P)$ can be written as
\begin{multline}\label{HS}
\varphi(h^2 P)=\sum_{k=0}^{N}h^k
a_{k}(x,hD)+
\\
+\frac{h^{N+1}}{\pi}\int_{\C}\bar{\partial}\widetilde{\varphi}(z)
(h^2 P-z)^{-1}r_{N+1}(x,hD,z,h)dL(z)\, 
\end{multline}
with
$$
a_{0}(x,\xi)=\varphi(p_0(x,\xi))
$$
and, for $k\geq 1$,
$$
a_{k}(x,\xi)=\sum_{j=1}^{2k-1}
\frac{(-1)^{j}}{j!} d_{j,k}(x,\xi)\varphi^{(j)}(p_0(x,\xi))\, .
$$
We now state a bound for the action of the map $a_{k}(x,hD)$ on Lebesgue spaces.
\begin{lemme}\label{oh}
For $1\leq q\leq r\leq\infty$,
\begin{equation*}
\|a_{k}(x,hD)\|_{L^q(\R^d)\rightarrow L^r(\R^d)}\leq C_{qrk}\,
h^{d(\frac{1}{r}-\frac{1}{q})}\, .
\end{equation*}
\end{lemme}
\begin{proof}
Since $a_k(x,\xi)$ is smooth and compactly supported with respect to $\xi$,
Lemma~\ref{oh} is a direct consequence of Proposition~\ref{sympa}.
\end{proof}
We next state a bound for the second term in the right hand-side of
(\ref{HS}).
\begin{lemme}\label{remainder}
Let $n_0>d/2$ be an integer. 
Set
$$
R_{N}:=\frac{h^{N+1}}{\pi}\int_{\C}\bar{\partial}\widetilde{\varphi}(z)
(h^2 P-z)^{-1}r_{N+1}(x,hD,z,h)dL(z)\, .
$$
Then for every $N>d+2n_0-3$, every $1\leq q\leq r\leq\infty$, there exists
$C>0$ such that for every $u\in
L^r(\R^d)$, every $h\in]0,1]$,
\begin{equation}\label{parvo}
\|R_{N}u\|_{L^q(\R^d)}\leq
Ch^{N+1-2n_0+d(\frac{1}{r}-\frac{1}{q})}\|u\|_{L^r(\R^d)}\,.
\end{equation}
Moreover for every $s\geq 0$ there exists $N_s$ such that for every $N\geq
N_s$, there exists $C>0$ such that for every $h\in]0,1]$, every $u\in H^{-s}(\R^d)$,
\begin{equation}\label{vtoro}
\|h^{N+1}R_{N}u\|_{H^s(\R^d)}\leq Ch^{N-N_s}\|u\|_{H^{-s}(\R^d)}\, .
\end{equation}
\end{lemme}
\begin{proof}
Define $\widetilde{r}_{N+1}(x,hD,z,h)$ by setting
$$
(h^2 P-z)^{-1}r_{N+1}(x,hD,z,h)=
(h^2 P-z)^{-n_0}\widetilde{r}_{N+1}(x,hD,z,h)\, ,
$$
with a fixed $n_0>d/2$.
Then $\widetilde{r}_{N+1}$ satisfies
$$
\big|
\partial_{x}^{\alpha}\partial_{\xi}^{\beta}\widetilde{r}_{N+1}(x,\xi,z,h)
\big|
\leq
C_{\alpha\beta}\langle\xi\rangle^{-N-1+2(n_0-1)-|\beta|}
\Big(\frac{\langle z\rangle}{|\Im\, z|}\Big)^{n(\alpha,\beta)}\, .
$$
If $N+1-2(n_0-1)>d$, i.e. $N>d+2n_0-3$, then Proposition~\ref{sympa} implies
that there exists $n_1\in\N$ such that
$$
\|\widetilde{r}_{N+1}(x,hD,z,h)\|_{L^q(\R^d)\rightarrow L^r(\R^d)}\leq C\,
h^{d(\frac{1}{r}-\frac{1}{q})}\,
\Big(\frac{\langle z\rangle}{|\Im\, z|}\Big)^{n_1}\, .
$$
On the other hand, we can write
$$
(h^2P-z)^{-n_0}=h^{-2n_0}(P-h^{-2}z)^{-n_0}
$$
thus Proposition~\ref{ppp} shows that 
there exists $n\in\N$ such that
for every $r\in[1,\infty]$,
\begin{equation*}
\|(h^2P-z)^{-n_0}\|_{L^r(\R^d)\rightarrow L^r(\R^d)}\leq
C_{r} h^{-2n_0}\Big(\frac{\langle h^{-2}z\rangle}{|\Im\, h^{-2} z|}\Big)^{n}
\leq 
\tilde{C}_{r}
h^{-2n_0}\Big(\frac{\langle z\rangle}{|\Im\, z|}\Big)^{n},
\end{equation*}
where we used that for every $h\in]0,1]$, every 
$z\in \C\backslash[0,\infty[$,
$$
\frac{\langle h^{-2}z\rangle}{|\Im\, h^{-2} z|}\leq C
\frac{\langle z\rangle}{|\Im\, z|}\,.
$$
The  proof of (\ref{parvo}) is completed by taking $\Lambda>n_1+n$ in (\ref{anal}).
Let us next prove (\ref{vtoro}). 
It suffices to prove the result for $s$ an even integer. It is sufficient to study the
action on $L^2(\R^d)$ of the map
$$
\int_{\C}\bar{\partial}\widetilde{\varphi}(z)
(h^2 P-z)^{-1}
(1+h^{-2}h^{2}P)^{s/2}
r_{N+1}(x,hD,z,h)
(1-h^{-2}h^{2}\Delta)^{s/2}
dL(z),
$$
where $\Delta$ is the flat Laplacian on $\R^d$. We can then write
$$
(1+h^{-2}h^{2}P)^{s/2}
r_{N+1}(x,hD,z,h)
(1-h^{-2}h^{2}\Delta)^{s/2}
=h^{-2s}\widetilde{r}_{N,s}(x,hD,z,h),
$$
where 
$\widetilde{r}_{N,s}(\cdot,\cdot,z,h)\in S^{-N-1+2s}$ with semi-norms 
{\it uniformly} bounded with respect to $h\in]0,1]$ by a 
by a polynomial of (\ref{frac}). Therefore, for $N+1-2s>d$, we can apply 
(\ref{im}) and Proposition~\ref{sympa} to conclude the proof of (\ref{vtoro}).
This completes the proof of  Lemma~\ref{remainder}.
\end{proof}
Combining Lemma~\ref{oh} and  Lemma~\ref{remainder} completes the proof of
Proposition~\ref{cf}.
\end{proof}
We will now give several consequences of Proposition~\ref{cf} that we will use
in the sequel. 
We first quote the following proposition which is a direct consequence of
Proposition~\ref{sympa} and Proposition~\ref{cf}.
\begin{proposition}\label{phi}
Let $\varphi\in C_{0}^{\infty}(\R)$.
Then for every $h\in]0,1]$, every $1\leq r\leq q\leq\infty$,
$$
\|\varphi(h^2 P)\|_{L^q(\R^d)\rightarrow L^r(\R^d)}\leq C_{qr}\,
h^{d(\frac{1}{r}-\frac{1}{q})}\,
$$
\end{proposition}
Next, we state a consequence of the Littlewood-Paley theory in terms of
$\varphi(h^2 P)$.
Consider a Littlewood-Paley partition of the
identity
$$
{\rm Id}\,=
\varphi_{1}(P)+\sum_{h^{-1}\,:\,\,{\rm dyadic}}\,\varphi(h^2P),
$$
where $\varphi_{1}\in C_{0}^{\infty}(\R)$, $\varphi\in
C_{0}^{\infty}(\R\backslash\{0\})$ and
``$h^{-1}\,:\,\,{\rm dyadic}$'' means that in the sum $h^{-1}$ takes all positive powers of
$2$ as values. 
The existence of such a partition is standard (see e.g. \cite{AG}).
We then have the following statement.
\begin{proposition}\label{LP}
Let $T>0$. Then for every $u\in C([0,T];{\mathcal S}(\R^d))$, every
$p\in[2,\infty]$, every $q\in [2,+\infty[$,
$$
\|u\|_{L^{p}_{T}L^q}\leq
C\|u\|_{L^{\infty}_{T}L^2}+C\big( \sum_{h^{-1}\,:\,\,{\rm dyadic}}\,
\|\varphi(h^2 P) u\|^{2}_{L^{p}_{T}L^q}\big)^{\frac{1}{2}}\, .
$$
\end{proposition}
\begin{proof}
Proceeding as in \cite[Corollary~2.3]{BGT1}, we obtain that
$$
\|u\|_{L^p_{T}L^q}\leq C\|u\|_{L^{\infty}_{T}L^2}+
C \Big\|
\big( \sum_{h^{-1}\,:\,\,{\rm dyadic}}\,
\|\varphi(h^2 P) u\|^{2}_{L^q}\big)^{\frac{1}{2}}
\Big\|_{L^p[0,T]}
\, .
$$
Since $p\geq 2$, the Minkowski inequality completes the proof of Proposition~\ref{LP}.
\end{proof}
In the proof of Theorem~\ref{thm2}, we will make use of the following statement.
\begin{proposition}\label{inverse}
Let $\varphi\in C_{0}^{\infty}(\R\backslash\{0\})$.
Then for every $s\in\R$ there exists $C_s$ such that for every 
$w\in {\mathcal S}(\R^d)$, every $h\in]0,1]$,
\begin{equation}\label{koko}
\|\varphi(h^2 P)w\|_{L^2(\R^d)}\leq C_{s}\,h^s\|w\|_{H^s(\R^d)}\, .
\end{equation}
\end{proposition}
\begin{proof}
Using (\ref{lesno}) and the spectral theorem for $P$, we can write
\begin{eqnarray*}
\|\varphi(h^2 P)w\|_{L^2} & = & \|\varphi(h^2 P)(1+P)^{-s/2}(1+P)^{s/2}w\|_{L^2}
\\
& \leq & C \|\varphi(h^2 P)(1+P)^{-s/2}\|_{L^2\rightarrow L^2}\|w\|_{H^s}
\\
& \leq & 
C\sup_{\lambda\sim h^{-2}}(1+\lambda)^{-s/2}\, \|w\|_{H^s}
\leq Ch^s\|w\|_{H^s}\,.
\end{eqnarray*}
This completes the proof of Proposition~\ref{inverse}.
\end{proof}
In applications to nonlinear problems, one may also need to use $L^p$ versions
of Proposition~\ref{inverse}.
Here is a precise statement.
\begin{proposition}\label{inverse-bis}
Let $\varphi\in C_{0}^{\infty}(\R\backslash\{0\})$ and $p\in[1,+\infty]$.
Then there exists $C$ such that for every 
$w\in {\mathcal S}(\R^d)$, every $h\in]0,1]$,
\begin{equation*}
\|\varphi(h^2 P)w\|_{L^p(\R^d)}\leq C\,h\|w\|_{W^{1,p}(\R^d)}\, .
\end{equation*}
\end{proposition}
\begin{proof}
In view of Proposition~\ref{cf}, it suffices to establish the bound
\begin{equation*}
\|a(x,hD)w\|_{L^p(\R^d)}\leq C\,h\|w\|_{W^{1,p}(\R^d)}\, ,
\end{equation*}
where $a(x,\xi)$ is satisfying (\ref{decay}) and (\ref{zero}). Let $\psi\in
C_{0}^{\infty}(\R^d)$ be such that
$$
a(x,\xi)\psi(\xi)=a(x,\xi)\, .
$$
Thanks to Proposition~\ref{sympa}, the map $a(x,hD)$ is bounded on $L^p(\R^d)$
and thus it is sufficient to prove that
\begin{equation*}
\|\psi(hD)w\|_{L^p(\R^d)}\leq C\,h\|w\|_{W^{1,p}(\R^d)}
\end{equation*}
which is a well-known fact (see e.g. \cite[Chapitre~2]{AG}, \cite{Chemin}).
This completes the proof of Proposition~\ref{inverse-bis}.
\end{proof}
\begin{remarque}
Proposition~\ref{inverse-bis} will be important in the proof of
Theorem~\ref{thm5} below. A similar bound in the context of a compact manifold
was used in \cite{BGT1}. The additional point here is again the $L^p$
boundedness of the remainder in the pseudo differential expansion of
$\varphi(h^2 P)$ established in Proposition~\ref{cf}.
\end{remarque}
\section{The Isozaki-Kitada parametrix}
This section is devoted to the construction of Isozaki-Kitada.
We only give the details for those arguments which are not written explicitly in the
papers of the reference list. The reader interested in having all the details for
the proofs of the statements in this section can consult \cite[Section~4]{Robert2},
\cite[Appendice]{bouclet},\cite[Appendix]{boucletBIS}.
The reader may of course wish to consult the original paper by Isozaki-Kitada \cite{IK}
which is nevertheless only written for potential perturbations and not in the
semi-classical regime.

In this section $g$ is a metric on $\R^d$ satisfying (\ref{H1}) and (\ref{H1bis}).
For $J \Subset ]0,+\infty[$ an open interval, $R>0$ and $\sigma\in]-1,1[$,
we consider the outgoing and incoming zones $\Gamma^{+}(R,J,\sigma)$ and
$\Gamma^{-}(R,J,\sigma)$,
defined by
$$
\Gamma^{\pm}(R,J,\sigma)=
\big\{
(x,\xi)\in\R^{2d}\,:\, |x|> R,\,\, |\xi|^{2}\in J,\,\,
\pm\frac{\langle x,\xi\rangle}{|x||\xi|}>-\sigma
\big\}\,.
$$
The next statement is proved in Robert \cite[prop. 4.1]{Robert2} (see also \cite{IK}).
\begin{proposition}\label{HJIK}
For every interval $J \Subset ]0,+\infty[$, every $\sigma\in]-1,1[$ there
exist a large number $R$ and phase functions 
$S_{\pm}\in C^{\infty}(\R^{2d};\R)$ satisfying
$$
g(x)(\nabla_{x}S_{\pm},\nabla_{x}S_{\pm})=|\xi|^{2},\quad
(x,\xi)\in \Gamma^{\pm}(R,J,\sigma)
$$
and, for every $(\alpha,\beta)\in\N^{2d}$, there exists $C_{\alpha\beta}$ such
that for every $(x,\xi)\in\R^{2d}$,
$$
\big|
\partial_{x}^{\alpha}\partial_{\xi}^{\beta}
\big(
S_{\pm}(x,\xi)-\langle x,\xi\rangle
\big)
\big|
\leq C_{\alpha\beta}\langle x\rangle^{1-\rho-|\alpha|}\, .
$$
\end{proposition}
We next state an easy consequence of Proposition~\ref{HJIK}.
\begin{proposition}\label{HJIKbis}
For every interval $J \Subset ]0,+\infty[$, every $\sigma\in]-1,1[$ there
exist a large number $\tilde{R}$ and two families of phase functions 
$$
(S_{\pm,R})_{R\geq \tilde{R}}\in C^{\infty}(\R^{2d};\R)
$$
such that
\begin{equation}\label{cin}
g(x)(\nabla_{x}S_{\pm,R},\nabla_{x}S_{\pm,R})=|\xi|^{2},\quad
(x,\xi)\in \Gamma^{\pm}(R,J,\sigma)
\end{equation}
and, such that 
for every $(\alpha,\beta)\in\N^{2d}$, there exists $C_{\alpha\beta}$ such
that for every $(x,\xi)\in\R^{2d}$, every $R\geq 2\tilde{R}$,
$$
\big|
\partial_{x}^{\alpha}\partial_{\xi}^{\beta}
\big(
S_{\pm,R}(x,\xi)-\langle x,\xi\rangle
\big)
\big|
\leq C_{\alpha\beta}
\min\Big(
\langle x\rangle^{1-\rho-|\alpha|},
R^{1-\rho-|\alpha|}
\Big).
$$
\end{proposition}
\begin{proof}
Let $J_0$ be an open interval such that $J\Subset J_0\Subset ]0,+\infty[$.
Define a smooth function $\theta_{J,J_0}$ such that 
$$
\theta_{J,J_0}(x)=
\left\{
\begin{array}{ll}
1 & \quad {\rm when}\quad x\in J,
\\
0 &\quad {\rm when}\quad x\notin J_0.
\end{array}
\right.
$$
Let us next fix $\sigma_0$ such that $\sigma_0\in ]\sigma,1[$.
Let $\kappa_{\sigma,\sigma_0}$ be a monotone smooth function such that
$$
\kappa_{\sigma,\sigma_0}(x)=
\left\{
\begin{array}{ll}
1 & \quad {\rm when}\quad x\geq -\sigma,
\\
0 &\quad {\rm when}\quad x\leq -\sigma_0.
\end{array}
\right.
$$
Next, let $\chi\in C^{\infty}$ be monotone and such that 
$$
\chi(x)=
\left\{
\begin{array}{ll}
1 & \quad {\rm when}\quad x\geq 1,
\\
0 &\quad {\rm when}\quad x\leq 1/2.
\end{array}
\right.
$$
For $R\geq 1$, we define $\chi_{R}(x)$ as 
$$
\chi_{R}(x)=\chi\Big(\frac{|x|}{R}\Big)\, .
$$
Then for $\tilde{R}\geq 1$ and $R\geq 2\tilde{R}$ the function $\psi_{R}(x,\xi)$ defined as
$$
\psi_{R}^{\pm}(x,\xi)=\chi_{R}(x)\theta_{J,J_0}(|\xi|^2)
\kappa_{\sigma,\sigma_0}\Big(\pm \frac{\langle x,\xi\rangle}{|x||\xi|}\Big)
$$
satisfies the properties
$$
{\rm supp}\,\,\psi_{R}^{\pm}\subset
\Gamma^{\pm}(R/2,J_0,\sigma_0)\subset\Gamma^{\pm}(\tilde{R},J_0,\sigma_0)
$$
and
$
\psi_{R}^{\pm}(x,\xi)=1
$
for $(x,\xi)\in \Gamma^{\pm}(R,J,\sigma)$.
The functions $\psi_{R}^{\pm}(x,\xi)$ also enjoy the bounds
\begin{equation}\label{trick}
\big|
\partial_{x}^{\alpha}\partial_{\xi}^{\beta}\psi_{R}^{\pm}(x,\xi)
\big|\leq
C_{\alpha\beta N}
\min\Big(
\langle x\rangle^{-|\alpha|},
R^{-|\alpha|}
\Big)
\langle\xi\rangle^{-N}
\end{equation}
with a constant $C_{\alpha\beta N}$ {\it uniform} with respect to $R\geq \tilde{R}$.

For $\tilde{R}\gg 1$, let us denote by $\tilde{S}_{\pm,\tilde{R}}(x,\xi)$ the phase function given
by Proposition~\ref{HJIK} associated to $J_0$, $\sigma_0$ and $\tilde{R}$.
Then, by invoking (\ref{trick}) 
and the fact that the derivatives of $\psi_{R}^{\pm}$ with respect to $x$ are
supported in a set $\{x \gtrsim  R\}$,
we observe that the phase functions $S_{\pm,R}$
defined as
$$
S_{\pm,R}(x,\xi)=
\psi_{R}^{\pm}(x,\xi)\tilde{S}_{\pm,\tilde{R}}(x,\xi)+
(1-\psi_{R}^{\pm}(x,\xi))\langle x,\xi\rangle
$$
satisfies the claimed properties.
This completes the proof of Proposition~\ref{HJIKbis}.
\end{proof}
For a given real number $\mu$, we denote by $S(\mu,-\infty)$ the set of
smooth functions $a(x,\xi)$ on $\R^{2d}$ such that for every $N\in \N$, every
$(\alpha,\beta)\in\N^{2d}$ there exists a constant $C_{N\alpha\beta}$ such that for every
$(x,\xi)\in\R^{2d}$,
$$
\big|\partial_{x}^{\alpha}\partial_{\xi}^{\beta}a(x,\xi)\big|
\leq 
C_{N\alpha\beta}\langle x\rangle^{\mu-|\alpha|}\langle\xi\rangle^{-N}\,.
$$
We equip $S(\mu,-\infty)$ with the natural Fr\'echet space topology.
The next proposition is devoted to the semi-classical Isozaki-Kitada parametrix.
\begin{proposition}\label{IKpar}
Let us fix an open interval $J\Subset ]0,+\infty[$ and $\sigma\in]-1,1[$. Consider
the open intervals $J_1$ and $J_2$ so that
$$
J\Subset J_1\Subset J_2 \Subset ]0,+\infty[
$$
and real numbers
$$
\sigma<\sigma_1<\sigma_2<1\,.
$$
Then there exists $R_0\gg 1$ such that for every $N\in\N \setminus \{ 0 \}   $, every $k\in\N$,
every $R\geq R_0$, every $\chi_{\pm}\in S(-k,-\infty)$ supported in
$\Gamma^{\pm}(R,J,\sigma)$ we can find:
\begin{itemize}
\item
a sequence 
$$
a_j^{\pm}\in S(-j,-\infty),\quad j=0,1,\dots, N
$$
of smooth functions on $\R^{2d}$
supported in $\Gamma^{\pm}(R^{1/3},J_2,\sigma_2)$, 
\item
a sequence 
$$
b_j^{\pm}\in S(-k-j,-\infty),\quad j=0,1,\dots, N
$$
of smooth functions on $\R^{2d}$
supported in $\Gamma^{\pm}(R^{1/2},J_1,\sigma_1)$
\end{itemize}
such that for every $h\in]0,1]$, every $\pm t\geq 0$,
\begin{eqnarray*}
e^{-i\frac{t}{h}h^2P}\chi_\pm(x,hD)
& = &
J_{S_{\pm,R^{1/4}}}
\Big(
\sum_{j=0}^{N}h^{j}a_{j}^{\pm}
\Big)e^{-i\frac{t}{h}h^2P_{0}}
J_{S_{\pm,R^{1/4}}}
\Big(
\sum_{j=0}^{N}h^{j}b_{j}^{\pm}
\Big)^{\star}
\\
& &
+
h^{N+1}R^{\pm}_{N}(t,h),
\end{eqnarray*}
where the phase functions $S_{\pm,R^{1/4}}$ are defined (with $R^{1/4}$ instead of $R$) in
Proposition~\ref{HJIKbis}, $P_0=-\Delta$ denotes the flat Laplacian on $\R^d$
and the maps
$J_{S_{\pm,R^{1/4}}}(q)$ are   defined by
$$
J_{S_{\pm,R^{1/4}}}(q)u(x)=
(2\pi h)^{-d}\int\int
e^{\frac{i}{h}(S_{\pm,R^{1/4}}(x,\xi)-\langle y,\xi\rangle)}\,
q(x,\xi)u(y)dyd\xi\, .
$$
Moreover, for every $T>0$, every $N\geq 1$ and every positive
integer $s$, 
there exists $C$ such that for
every $h\in]0,1]$, every $\pm t\in[0,Th^{-1}]$,
the remainder $R^{\pm}_{N}(t,h)$ satisfies 
\begin{equation}\label{vac}
\Big\|(P+1)^{s}R^{\pm}_{N}(t,h)(P+1)^{s}\Big\|_{L^2(\R^d)\rightarrow L^2(\R^d)}
\leq C\, h^{N-4s-2}\,.
\end{equation}
\end{proposition}
\begin{proof}
For the precise construction of  $a_j^{\pm} $ and $ b_j^{\pm} $ 
we refer to \cite[Section~4]{Robert2},\cite{bouclet}.
However, for the convenience of the reader, we recall the main
lines of the method in the outgoing case (the incoming one being
similar). 
We first choose $ J_3 , \sigma_3 $ such that 
$ J_2\Subset J_3 $, $ \sigma_2 < \sigma_3 < 1 $ and then choose
$S_{+,R^{1/4}}\equiv S_+ $ (for shortness) as in Proposition~\ref{HJIKbis} solving
the Hamilton-Jacobi equation (\ref{cin}) on 
$ \Gamma^+ (R^{1/4},J_3,\sigma_3)$. We then look for a symbol
$$ 
a^+ = a_0^+ + h a_1^+ + \cdots + h^N a_N^+
$$ 
supported in $\Gamma^+ (R^{1/3},J_2,\sigma_2) $ such that
$$  
J_{S_+} (c_N(h)) :=  (h^2 P) J_{S_+}(a^+) - J_{S_+}(a^+) (h^2 P_0)   
$$
has a ``small contribution'' (see (\ref{duhamel}) below).
This leads to a system of equations for
$ a_0^+ , \cdots , a_N^+ $ that take the form of rather standard
(time independent) transport equations, but only in the region
where the Hamilton-Jacobi equation (\ref{cin}) is satisfied.
For $ R $ large enough, we can solve these equations in a
neighborhood of $ \Gamma^+ (R^{1/3},J_2,\sigma_2) $ and by cutting
off these solutions by a function supported in 
$ \Gamma^+(R^{1/3},J_2,\sigma_2) $ which equals $1$ near 
$ \Gamma^+(R^{5/12},J_1,\tilde{\sigma}_1) $ (notice that $ 1 /3 < 5/12 < 1/2  $) with
$  \sigma_1  < \tilde{\sigma}_1 < \sigma_2  $,
we can build $ a_0^+, \cdots , a_N^+ $
so that
\begin{eqnarray}
c_N (h) = h^{N+1} \tilde{r}_{N+1}(h) + \tilde{c}_N (h)
\end{eqnarray}
with $ (\tilde{r}_{N+1}(h))_{0<h \leq 1} $ in a bounded set of 
$ S(-N,-\infty) $ and $ (\tilde{c}_N (h))_{0<h\leq 1} $ in a bounded
set of $ S (0, -\infty) $, supported in $ \Gamma^+
(R^{1/3},J_2,\sigma_2) $
 and such that
\begin{eqnarray}
 \tilde{c}_N (h) \equiv 0 \qquad \mbox{near} \ \ \Gamma^+
(R^{5/12},J_1, \tilde{\sigma}_1) . \label{nonstationnaire}
\end{eqnarray}
Note that the symbol $ c_N (h) $ is a priori not small, because of the term
$ \tilde{c}_N (h)$. However, using (\ref{nonstationnaire}),
we shall see afterwards that the contribution of $ J_{S_+}(
\tilde{c}_N (h) ) $ to the final remainder term of the parametrix
is harmless, once multiplied from the right by another FIO with nice support
properties. Next, remarking that we can choose $ a_0^+ $ non
vanishing near $ \Gamma^+ (R^{1/2},J_1,\sigma_1) $, we can solve
another family of (algebraic) equations for $ b_0^+ , \cdots ,
b_N^+ $ such that 
$$ 
b^+ := b_0^+ + \cdots + h^N b_N^+ 
$$ 
is
supported in $ \Gamma^+ (R^{1/2},J_1,\sigma_1) $ and satisfies
$$ 
J_{S_+} (a^+) J_{S^+} (b^+)^{\star} = \chi_+ (x,hD) + h^{N+1} r_N
(x,hD,h) 
$$
with $ (r_N (h))_{0 < h \leq 1} $ in a bounded set of $ S(-N,-\infty) $.
More precisely, the standard composition rule for the computation of the
symbol of  $ J_{S_+}(a^+) J_{S^+}(b^+)^{\star}  $ (see \cite{Robert}) show that the
above equation leads to a triangular system with unknown  $ b_0^+ , \ldots
,b_N^+  $
and $ a_0^+  $ on the diagonal.

Combining this last equation and the fact that
\begin{multline}\label{duhamel}
e^{-i \frac{t}{h} h^2 P} J_{S_+} (a_+) - J_{S_+} (a_+) e^{-i
\frac{t}{h} h^2 P_0} 
\\
=
-
\frac{i}{h} \int_0^t e^{-i \frac{t - \tau}{h} h^2 P} J_{S_+}
(c_N(h)) e^{-i \frac{\tau}{h} h^2 P_0}  d \tau   , 
\end{multline}
we can then represent the remainder of the Isozaki-Kitada parametrix as
$$
h^{N+1}R^{+}_{N}(t,h) =I+II+III,
$$
where
\begin{eqnarray*}
I & = & -h^{N+1}e^{-i\frac{t}{h}h^2P}r_{N+1}(x,hD,h),
\\
II & = & -ih^{N} \int_{0}^{t}e^{-i\frac{t-\tau}{h}h^2P}
J_{S_{+}}(\tilde{r}_{N+1}(h))e^{-i\frac{\tau}{h}h^2P_0}d\tau \,\circ
\, J_{S_{+}} \big( b^+  \big)^{\star}
\\
III & = &
-\frac{i}{h}\int_{0}^{t}e^{-i\frac{t-\tau}{h}h^2P}K_{h}(\tau)d\tau,
\end{eqnarray*}
with  $K_{h}(\tau)$  an operator with kernel $ {\mathcal
K}(x,y,h,\tau) $ satisfying, for all $ M \geq 0 $ and $
\alpha,\beta \in {\mathbb N}^d $,
\begin{eqnarray}
\big|\partial_{x}^{\alpha}\partial_{y}^{\beta} {\mathcal
K}(x,y,h,\tau) \big|\leq C_{\alpha\beta
M}h^{M}(1+|x|+|y|+\tau)^{-M} , \qquad \tau \geq 0 . \label{phasenonstationnaire}
\end{eqnarray}
Here $ K_h (\tau) = J_{S_+} (\tilde{c}_N(h))  e^{i \frac{\tau}{ h}
h^2 P_0} J_{S_+} (b^+)^{\star} $ contains the contribution of $
\tilde{c}_N(h) $. The above estimate follows from a non stationary
phase argument (see \cite[Proposition~2.4.7]{bouclet}), by exploiting the
support properties of $ \tilde{c}_N (h) $ and $ b^+ $. For the sake of
completeness, and to somehow prepare the reader to the stationary phase argument used in
the next section, we recall how to prove $ ( \ref{phasenonstationnaire} )
$. This proof can be found in \cite{IK} for potential perturbations in the non
semi-classical case ($h=1$). Here we reproduce the proof of \cite{bouclet}. 

The kernel of $K_{h}(\tau)$ at the point $(x,y)$ is given by an oscillatory
integral over a fixed compact set in the $\xi$ variable and the phase function
$$
ih^{-1}(\tau |\xi|^2 + S_+ (y,\xi) - S_+(x,\xi)).
$$
More precisely the integration is over those $\xi$ such that
\begin{equation}\label{conditionsupport} 
(x,\xi) \in   \Gamma^+ (R^{1/3},J_2,\sigma_2)\setminus
\Gamma^+ (R^{5/12},J_1 , \tilde{\sigma}_1)  \ \ \mbox{and} \ \   (y,\xi) \in  \Gamma^+
(R^{1/2},J_1,\sigma_1) .
\end{equation}
Estimate $(\ref{phasenonstationnaire})$ 
follows from standard integrations by parts using the following lemma.
\begin{lemme}\label{ddds}
There exist $c > 0$ and $R_0>1$ such that for every $ \tau \geq 0  $, every
$x$, $y$, $\xi$ satisfying (\ref{conditionsupport}), every $R\geq R_0$,
\begin{eqnarray}
 \left|  \nabla_{\xi} ( \tau |\xi|^2 + S_+ (y,\xi) - S_+(x,\xi) )  \right| \geq
c (1 + \tau + |x| + |y|).
\label{minorationphase}
\end{eqnarray}
\end{lemme}
\begin{proof} 
For simplicity, we use the notation $ \mbox{cos}(x,\xi) = \langle x , \xi \rangle / |x||\xi|  $.  We will also use the following statements that are easily checked:
\begin{eqnarray}
\cos(\eta,\eta^{\prime}) \geq - \sigma>-1 \Rightarrow |\eta+\eta^{\prime}| \geq
(1-|\sigma|)^{1/2}(|\eta|^2 + |\eta^{\prime}|^2)^{1/2}, \nonumber \\
|\eta-\eta^{\prime}| \leq \epsilon |\eta| \Rightarrow |\cos ( \eta^{\prime \prime},\eta) -
\cos(\eta^{\prime \prime},\eta^{\prime})| \leq 2 \epsilon, \ \ \forall
\eta^{\prime \prime}
\in {\mathbb R}^d \setminus 0 . \nonumber
\end{eqnarray}
One then remark that if $(\ref{conditionsupport})$ holds then either $   |x| \leq R^{5/12} $
or $ - \tilde{\sigma}_1 \geq \cos (x,\xi) > - \sigma_2 $. 

Assume first that $ |x|\leq R^{5/12}  $. Fix $ \epsilon  $ such that 
$ \sigma_2 - \sigma_1 > 2 \epsilon  $.  
By choosing $ R $ large enough, we have $ |\nabla_{\xi} S_+
(y,\xi)-y| \leq \epsilon |y|  $ and thus 
$  \cos (\nabla_{\xi} S_+ (y,\xi) , \xi  ) > - |\sigma_2|  $. This implies that
$$ |\nabla_{\xi} (\tau |\xi|^2 + S_+ (y,\xi))  | \geq (1-|\sigma_2|)^{1/2} (4
\tau^2 |\xi|^2 + |\nabla_{\xi} S_+ (y,\xi)|^2)^{1/2}  $$
from which $ (\ref{minorationphase})  $ follows easily since 
$ | \nabla_{\xi} S_+ (x,\xi)  | \lesssim R^{5/12}  $ and $ | \nabla_{\xi} S_+
(y,\xi) | \gtrsim R^{1/2}  $.

We now  assume that $ - \tilde{\sigma}_1 \geq \cos (x,\xi) > - \sigma_2    $.
It suffices to show that, if $R$ is large enough, there exists  $ \sigma \in (-1,1)  $ such
that
\begin{equation}
 \cos \left( \nabla_{\xi} (S_+(y,\xi)-S_+(x,\xi)) , \xi  \right) \geq \sigma , 
\label{premier}
\end{equation}
\begin{equation}
\cos \left( -  \nabla_{\xi} S_+ (x,\xi) , \nabla_{\xi} S_+ (y,\xi)  \right) \geq
\sigma 
\label{deuxieme}
\end{equation}
provided   $ (\ref{conditionsupport}) $ holds. 
Indeed, the estimate $(\ref{premier})$ implies that the left hand side of $
(\ref{minorationphase})  $ is bounded from below by $ c (4 \tau|\xi|^2 + |\nabla_{\xi}(S_+(x,\xi)-S_+(y,\xi))| )  $
and then $(\ref{deuxieme})$ implies that $ |\nabla_{\xi}(S_+(x,\xi)-S_+(y,\xi))|
\geq c ( |x|+|y| )$. 

It remains to prove $(\ref{premier})$ and $(\ref{deuxieme})$.
Let  us choose $ \epsilon > 0  $ such that $ \tilde{\sigma}_1 - \sigma_1 > 4 \epsilon$.
Observe that $ \sigma_1 + 2 \epsilon \in (-1,1)  $. By choosing $ R  $ large enough, we
have
\begin{equation} 
\cos(\nabla_{\xi} S_+ (y,\xi),\xi) \geq - \sigma_1 - 2 \epsilon  \qquad
\mbox{and} \qquad  2 \epsilon - \tilde{\sigma}_1  > \cos
(\nabla_{\xi} S_+ (x,\xi),\xi)  
\label{trigonometrie}
\end{equation}
Therefore
$
\nabla_{\xi} S_+ (y,\xi)\neq \nabla_{\xi} S_+ (x,\xi)  
$
and the left hand side of $ ( \ref{premier} ) $ reads
\begin{eqnarray}
\frac{|\nabla_{\xi}S_+ (y,\xi) | \cos (\nabla_{\xi}S_+(y,\xi),\xi) -
  |\nabla_{\xi}S_+(x,\xi)|
\cos (\nabla_{\xi}S_+(x,\xi),\xi) }{|\nabla_{\xi}S_+(y,\xi)- \nabla_{\xi}
S_+(x,\xi) |}  .
\nonumber
\end{eqnarray}
Now, using that for every $\alpha\in\R$, every $X,Y\in\R^d$, $X\neq Y$,
$$
\alpha\frac{|Y|-|X|}{|Y-X|}\geq -|\alpha|
$$
we deduce that $ (\ref{premier})  $ holds with $ \sigma = -| \sigma_1 + 2
\epsilon |  $. 
Finally, we see that $ (\ref{deuxieme}) $ must
hold for some possibly lower $ \sigma $ since otherwise 
we could find sequences $ (x_j), (y_j), (\xi_j)  $ satisfying
$ (\ref{conditionsupport})  $, $ - \tilde{\sigma}_1 \geq   \cos (x_j , \xi_j)
> - \sigma_2  $ and such that
$$ \lim_{j \rightarrow \infty} \frac{\nabla_{\xi} S_+ (x_j,\xi_j) }{
|\nabla_{\xi} S_+ (x_j,\xi_j)|} = \lim_{j \rightarrow \infty} 
\frac{\nabla_{\xi} S_+ (y_j,\xi_j)}{|\nabla_{\xi} S_+ (y_j,\xi_j)|}  $$
which is forbidden by $ ( \ref{trigonometrie} ) $ and the fact that $ -
\sigma_1 - 2 \epsilon > 2 \epsilon - \tilde{\sigma}_1  $.
This completes the proof of Lemma~\ref{ddds}.
\end{proof}
Therefore, we get the bounds
$$
\|(P+1)^{s}r_{N+1}(x,hD,h)(P+1)^{s}\|_{L^2\rightarrow L^2}\leq
Ch^{-4s} ,  
$$
and
$$
\|(P+1)^{s}K_{h}(\tau)(P+1)^{s}\||_{L^2\rightarrow L^2}\leq Ch^{M} ,
$$
for $\tau\geq 0$ and $h\in]0,1]$  
(in particular, we see {\it a posteriori} that the contribution of $ \tilde{c}_N (h) $ in the
remainder of the parametrix is $ {\mathcal O} (h^{\infty}) $).
Finally, using the $L^2$ boundedness of FIO (see e.g.
\cite{Robert}), we have that, for every $a\in S(0,-\infty)$ and
every $k\in {\mathbb N}$, there exists $C_{k,a}$ such that
$$
\|(h^2 P)^{k} J_{S_+ }(a)\|_{L^2\rightarrow L^2}\leq C_{k,a},
$$
and hence
\begin{multline*}
\Big\|(P+1)^{s}
J_{S_{+}}(\tilde{r}_{N+1}(h))e^{-i\frac{\tau}{h}h^2P} J_{S_{+}}
\big( b^+ \big)^{\star} (P+1)^{s} \Big\|_{L^2\rightarrow L^2}
\\
= \Big\|(P+1)^{s}
J_{S_{+}}(\tilde{r}_{N+1}(h))e^{-i\frac{\tau}{h}h^2P} \big(
(P+1)^{s} J_{S_{+}} (b^+) \big)^{\star} \Big\|_{L^2\rightarrow L^2}
\\
\leq Ch^{-4s}\,.
\end{multline*}
By integrating the corresponding estimates to $II$ and $III$
over an interval of size $T h^{-1} $ we get  the result (\ref{vac}).
\end{proof}
\begin{remarque}
Note that the control of the remainder is easier in our case than in
\cite{Robert2,bouclet} since we only need to integrate on $[0,Th^{-1}]$ 
in $\tau$ whereas in \cite{Robert2,bouclet} one has to integrate over
$\R^{+}$. Moreover, we do not use any non trapping assumption on the metric:
this is the main point in this paper.
\end{remarque}
\section{Strichartz estimates outside a large ball}
The goal of this section is to prove (\ref{est1}). The main point is to prove
the following statement.
\begin{proposition}\label{est1prelim}
Let $\varphi\in C_{0}^{\infty}(\R\setminus \{0\})$ and let $g$ be a metric on $\R^d$ satisfying (\ref{H1}) and (\ref{H1bis}). 
Then there exists $R>0$ such that  
for every $T>0$, every $(p,q)$ satisfying (\ref{adm}), every $\chi\in
C_{0}^{\infty}(\R^d)$, $\chi\equiv 1$ for $|x|<R$,
there exists  $C>0$ such that for every $f\in
L^2(\R^d)$, every $h\in]0,1]$,
\begin{equation}\label{ufr}
\|(1-\chi)e^{-itP}\varphi(h^2 P)
f\|_{L^{p}([-T,T];L^{q}(\R^d))}
\leq
C\|f\|_{L^2(\R^d)}\, .
\end{equation} 
\end{proposition}
\begin{remarque}\label{localized}
We could have spectrally localized $f$ in the right hand-side of (\ref{ufr}), i.e.
one could have replaced $\|f\|_{L^2(\R^d)}$ by $\|\varphi(h^2
P) f\|_{L^2(\R^d)}$. Indeed, let $\tilde{\varphi}\in
C_{0}^{\infty}(\R)$ be equal to one on the support of $\varphi$.
Then $\varphi=\varphi \tilde{\varphi}$ and we can apply (\ref{ufr})
with $\tilde{\varphi}$ instead of $\varphi$.
\end{remarque}
\begin{proof}
Recall that we denote by $P$ the self adjoint realization of $-\Delta_{g}$ on
$L^2(\R^d)$. 
For $\chi\in C_{0}^{\infty}(\R^d)$,
$\chi\equiv 1$ for $|x|<R$,
a partition of unity argument and Proposition~\ref{cf} allow to write
$$
(1-\chi)\varphi(h^2 P)=\sum_{k=0}^{N}h^{k}
\Big(\theta^{+}_{k}(x,hD)+\theta^{-}_{k}(x,hD)\Big)+h^{N+1}R_{N,\chi}(h)
$$
where $\theta^{\pm}_{k}\in S(-k,-\infty)$ and for some $J\Subset ]0,+\infty[$,
$\sigma_{\pm}\in  ]-1,1[$,
$$
{\rm supp}\,\, \theta_{k}^{\pm}\subset \Gamma^{\pm}(R,J,\sigma_{\pm})\, .
$$
More precisely, $\sigma_{\pm}$ and $J$ should be such that
$$
\big\{
(x,\xi)\,:\,
|x|>R,\,\,
p(x,\xi)\in {\rm supp}\,(\varphi)
\big\}
\subset
\Gamma^{+}(R,J,\sigma_{+})
\cup
\Gamma^{-}(R,J,\sigma_{-}).
$$
Let us notice that $\sigma_{\pm}$ can be taken both $1/2$.

Furthermore the remainder is such that for every $s\geq 0$ 
there exists $C>0$ such that for every $h\in ]0,1]$,
$$
\|(P+1)^{s/2}R_{N,\chi}(h)(P+1)^{s/2}\|_{L^2\rightarrow L^2}\leq Ch^{-2s}\,.
$$
In addition by the elementary properties of the $h$ pseudo differential
calculus (cf. e.g. \cite{Robert}), we can also write
$$
(1-\chi)\varphi(h^2 P)=\sum_{k=0}^{N}h^{k}
\Big(\chi^{+}_{k}(x,hD)^{\star}+\chi^{-}_{k}(x,hD)^{\star}\Big)+h^{N+1}\tilde{R}_{N,\chi}(h)\, ,
$$
where $\chi_{k}^{\pm}$ and $\tilde{R}_{N,\chi}(h)$ have similar properties to
$\theta_{k}^{\pm}$ and $R_{N,\chi}(h)$ respectively. Using the Sobolev
embedding, by taking $N$ large enough, we get the bound
$$
\big\|h^{N+1}\tilde{R}_{N,\chi}(h)e^{-itP}f\big\|_{L^p_{T}L^q}\leq C\|f\|_{L^2},
$$
provided $(p,q)$ is satisfying (\ref{adm}).
Therefore it suffices to prove the bound
\begin{equation}\label{etoile}
\big\|\chi^{\pm}_{k}(x,hD)^{\star}\,e^{-itP}f\big\|_{L^p_{T}L^q}\leq
C\|f\|_{L^2}\, .
\end{equation}
Since $\chi^{\pm}_{k}(x,hD)^{\star}\,e^{-itP}$ 
are clearly $L^2$ bounded, uniformly in $h$ and $t$, thanks
to the Keel-Tao theorem (see \cite{KT}, \cite[Proposition 2.8]{BGT1}), 
to get (\ref{etoile}), it suffices to prove the dispersive inequality
$$
\Big\|
\chi^{\pm}_{k}(x,hD)^{\star}\,e^{-itP}\,
e^{isP}\chi^{\pm}_{k}(x,hD)f
\Big\|_{L^{\infty}}
\leq 
\frac{C}{|t-s|^{d/2}}\|f\|_{L^1(\R^d)},\quad t,s\in [-T,T]\, ,
$$
uniformly with respect to $h$.
By the time rescaling $t\mapsto ht$, 
and by
defining the maps
$$
U^{\pm}_{k,h}(t)=\chi^{\pm}_{k}(x,hD)^{\star}\,e^{-ihtP},\quad t\in
[-h^{-1}T,h^{-1}T]\,,
$$
it suffices to prove the dispersive inequality
\begin{equation}\label{disp}
\big\|U^{\pm}_{k,h}(t)\, \big(U^{\pm}_{k,h}(s)\big)^{\star}\, f\big\|_{L^{\infty}(\R^d)}
\leq
\frac{C}{(h|t-s|)^{d/2}}\|f\|_{L^1(\R^d)},\quad t,s\in [-h^{-1}T,h^{-1}T]\,.
\end{equation}
Clearly
\begin{equation}\label{calcul}
U^{\pm}_{k,h}(t)\,
\big(U^{\pm}_{k,h}(s)\big)^{\star}=\chi_{k}^{\pm}(x,hD)^{\star}\,
e^{-ih(t-s)P}\chi_{k}^{\pm}(x,hD)\,.
\end{equation}
Denote by $K_{\pm}(t-s,x,y,h)$ (we do not explicit the dependence on $k$)
the kernel of (\ref{calcul}). In order to prove (\ref{disp}), it is sufficient
to show that there exists $C>0$ such that for every $h\in]0,1]$, every
$x,y\in\R^d$, every $t,s\in [-h^{-1}T,h^{-1}T]$, $t\neq s$, 
\begin{equation}\label{dispbis}
|K_{\pm}(t-s,x,y,h)|\leq C|h(t-s)|^{-d/2}\, .
\end{equation}
The next lemma contains the main trick in our analysis.
\begin{lemme}\label{tricky}
For $K_{+}$,
it suffices to prove (\ref{dispbis}) for $t-s>0$.
Similarly, for  $K_{-}$, it suffices to prove (\ref{dispbis}) for $t-s<0$.
\end{lemme}
\begin{proof}
We only consider $K_{+}$, the analysis for $K_{-}$ being similar.
Suppose that (\ref{dispbis}) holds true for $t-s>0$. Let 
$t,s\in  [-h^{-1}T,h^{-1}T]$ such that $t-s<0$.
Since
$$
U^{+}_{k,h}(t)\,
\big(U^{+}_{k,h}(s)\big)^{\star}=
\Big(
U^{+}_{k,h}(s)\,
\big(U^{+}_{k,h}(t)\big)^{\star}
\Big)^{\star}
$$
we obtain that
$$
K_{+}(t-s,x,y,h)=\overline{K_{+}(s-t,y,x,h)}\,.
$$
Since $s-t>0$, our assumption that (\ref{dispbis}) holds for positive values
of the first argument of $K_{+}$ implies that
$$
|K_{+}(t-s,x,y,h)|=
|\overline{K_{+}(s-t,y,x,h)}|\leq C(h(s-t))^{-d/2}
$$
which completes the proof of Lemma~\ref{tricky}.
\end{proof}
It is now clear that the proof of (\ref{dispbis}), and thus of 
Proposition~\ref{est1prelim}, will be finished, once we establish the
following lemma.
\begin{lemme}\label{nedelia}
There exists $R\gg 1$ and $C>0$ such that for every $h\in]0,1]$, every $\pm t\in
]0,h^{-1}T]$, every $x,y\in\R^d$,
$$
|K_{\pm}(t,x,y,h)|\leq C(\pm ht)^{-d/2}\,.
$$
\end{lemme}
\begin{proof}
As before, we only consider the case of $K_{+}$.
Denote by ${\mathcal K}_{+}(t,x,y,h)$ the kernel of
$\exp(-ithP)\chi_{k}^{+}(x,hD)$ and by $\tilde{{\mathcal K}}_{+}(x,y,h)$ the kernel of
$\chi_{k}^{+}(x,hD)^{\star}$. 
Clearly $\tilde{{\mathcal K}}_{+}(x,y,h)$ satisfies the assumptions of
the Schur lemma uniformly in $h$, and therefore by writing
$$
K_{+}(t,x,y,h)=\int_{\R^d}
\tilde{{\mathcal K}}_{+}(x,z,h)
{\mathcal K}_{+}(t,z,y,h)dz,
$$
we infer that
$$
|K_{+}(t,x,y,h)|\leq \sup_{z\in\R^d}|{\mathcal K}_{+}(t,z,y,h)|
\int_{\R^d}
|\tilde{{\mathcal K}}_{+}(x,z,h)|dz
\leq 
C\sup_{z\in\R^d}|{\mathcal K}_{+}(t,z,y,h)|\,.
$$
Therefore, it suffices to prove that there exists $C>0$ such that for every 
$h\in]0,1]$, every $t\in
]0,h^{-1}T]$, every $x,y\in\R^d$,
\begin{equation}\label{disptris}
|{\mathcal K}_{+}(t,x,y,h)|\leq C(ht)^{-d/2}\,.
\end{equation}
In order to prove (\ref{disptris}), we will of course use
Proposition~\ref{IKpar}.
With the notations of Proposition~\ref{IKpar}, since for $N\gg 1$, the map
$R_{N}(t,h)$ is bounded, uniformly with respect to $t,h$, from $H^{-s}$ to
$H^s$ with $s>d/2$, we deduce that its kernel is bounded uniformly with
respect to $h\in]0,1]$, $t\in[0,h^{-1}T]$. Therefore, in view of
Proposition~\ref{IKpar},
by expressing the kernel of 
$$
J_{S_{+,R^{1/4}}}
\big(a_{j_1}^{+}\big)
e^{-i\frac{t}{h}h^2P_{0}}
J_{S_{+,R^{1/4}}}
\big(b_{j_2}^{+}\big)^{\star},\quad 0\leq j_1,j_2\leq N,
$$
estimates (\ref{disptris}) will be established once we prove that 
there exists $C>0$ such that for every 
$h\in]0,1]$, every $t\in
]0,h^{-1}T]$, every $x,y\in\R^d$,
\begin{equation}\label{disp4}
(2\pi h)^{-d}\,
\Big|
\int_{\R^d}
e^{\frac{i}{h}\Phi_{+}(t,R,x,y,\xi)}a(x,\xi)\overline{b(y,\xi)}d\xi \Big|\leq C(ht)^{-d/2}\,,
\end{equation}
where the phase $\Phi_{+}$ is  defined as
$$
\Phi_{+}(t,R,x,y,\xi)=S_{+,R^{1/4}}(x,\xi)-S_{+,R^{1/4}}(y,\xi)-t|\xi|^{2}
$$
and $a(x,\xi)$, $b(x,\xi)$ are fixed smooth functions supported in
$$
\{(x,\xi)\in\R^{2d}\, :\, |x|\geq R\gg 1,\,\, 0<c\leq |\xi|\leq C\}
$$
for some positive constants $c$ and $C$.
In the proof of (\ref{disp4}), we will consider two different regimes for
$t$. If $t\in[0,h]$ then, using the support property with respect to $\xi$, 
the left hand-side of (\ref{disp4}) can be estimate by $Ch^{-d}$ which, for
$t\leq h$, is bounded by $C(ht)^{-d/2}$. Therefore, we can suppose that $t\geq
h$ in (\ref{disp4}). In this case we will take advantage of the rapid oscillations
of $\exp(ih^{-1}\Phi_{+})$. When $t\geq h$, the natural big parameter is
$th^{-1}$. We thus set
$$
\frac{1}{h}\Phi_{+}(t,R,x,y,\xi)=\frac{t}{h}\,\tilde{\Phi}_{+}(t,R,x,y,\xi),
$$
where
$$
\tilde{\Phi}_{+}(t,R,x,y,\xi)=\frac{S_{+,R^{1/4}}(x,\xi)-S_{+,R^{1/4}}(y,\xi)}{t}-|\xi|^{2}\,.
$$
We can write
$$
\tilde{\Phi}_{+}(t,R,x,y,\xi)=
\int_{0}^{1}
\langle\nabla_{x}\,
S_{+,R^{1/4}}(y+\alpha(x-y),\xi),x-y\rangle\frac{d\alpha}{t}-|\xi|^{2}\, .
$$
Therefore, thanks to the properties of the phase function $S_{+,R^{1/4}}$
displayed in Proposition~\ref{HJIKbis}, we obtain that
$$
\nabla_{\xi}\,
\tilde{\Phi}_{+}(t,R,x,y,\xi)=\frac{x-y}{t}-2\xi+Q(R,x,y,\xi)\cdot
\frac{x-y}{t},
$$
where $Q$ is a $d\times d$ matrix satisfying the bound
$$
|\partial_{\xi}^{\beta}Q(R,x,y,\xi)|\leq C_{\beta}R^{-\rho},
$$
for $R\gg 1$, $x,y\in\R^d$, $\xi$ on the support of $a(x,\xi)b(y,\xi)$.

Therefore, there exist $R_0\gg 1$ and $C_0>0$ such that for $R>R_0$ and
$\big|\frac{x-y}{t}\big|\geq C_0$, we have that, for $|\xi|\in [c,C]$,
\begin{equation}\label{Est1}
|\nabla_{\xi} \tilde{\Phi}_{+}(t,R,x,y,\xi)|\gtrsim \big|\frac{x-y}{t}\big|
\end{equation}
and 
\begin{equation}\label{Est2}
|\partial_{\xi}^{\beta}\tilde{\Phi}_{+}(t,R,x,y,\xi)|\leq
C_{\beta}\big|\frac{x-y}{t}\big|,\quad |\beta|\geq 2.
\end{equation}
Therefore using (\ref{Est1}), (\ref{Est2}) and integration by parts with
respect to $\xi$, we deduce that for every $\Lambda\geq 1$,
the left hand-side of (\ref{disp4}) is bounded by
$$
Ch^{-d}\, C_{\Lambda}(th^{-1})^{-\Lambda}\leq C(ht)^{-d/2}\, ,
$$
provided $\Lambda$ is taken bigger than $d/2$. 

We can therefore suppose that $t,x,y$ are such that $t\geq h$ and
\begin{equation}\label{boundimp}
\big|\frac{x-y}{t}\big|\leq C_0\,.
\end{equation}
In this case, we evaluate the left hand-side of (\ref{disp4}) by the
stationary phase. Under the condition (\ref{boundimp}), we can write
$$
\nabla_{\xi}^{2} \tilde{\Phi}_{+}(t,R,x,y,\xi)=-2{\rm Id}\,
+\nabla_{\xi}Q(R,x,y,\xi)\cdot \frac{x-y}{t}=
-2{\rm Id}\,+{\mathcal O}(R^{-\rho})\,.
$$
Therefore for $R\gg 1$, the map
$$
\xi\longmapsto \nabla_{\xi} \tilde{\Phi}_{+}(t,R,x,y,\xi)
$$
is a diffeomorphism from $\R^d$ to $\R^d$. In particular, 
for fixed $t\geq h$, $x$, $y$ satisfying (\ref{boundimp}), 
the phase $\tilde{\Phi}_{+}$ has a unique non degenerate critical point
$\xi_{cr}(t,R,x,y)$. Moreover, thanks to (\ref{boundimp}), for 
$|\beta|\geq 1$,
$$
\Big|
\partial_{\xi}^{\beta}\tilde{\Phi}_{+}(t,R,x,y,\xi_{cr}(t,R,x,y))
\Big|\leq C_{\beta}\, .
$$
We can therefore apply the stationary phase estimate to conclude that for
$t\geq h$ and $(t,x,y)$ satisfying (\ref{boundimp}), the left hand-side of
(\ref{disp4}) is bounded by 
$$
Ch^{-d}\, (th^{-1})^{-d/2}=C(ht)^{-d/2}\,.
$$
This completes the proof of Lemma~\ref{nedelia}.
\end{proof}
This completes the proof of Proposition~\ref{est1prelim}
\end{proof}
It is now clear that (\ref{est1}) will be proved, once we establish the
following statement.
\begin{proposition}\label{est1prelimbis}
Let $g$ be a metric on $\R^d$ satisfying (\ref{H1}) and (\ref{H1bis}). 
Then there exists $R>0$ such that  
for every $T>0$, every $(p,q)$ satisfying (\ref{adm}), every $\chi\in
C_{0}^{\infty}(\R^d)$, $\chi\equiv 1$ for $|x|<R$,
there exists  $C>0$ such that for every $f\in L^2(\R^d)$, 
\begin{equation*}
\|(1-\chi)e^{-itP} f\|_{L^{p}([-T,T];L^{q}(\R^d))}\leq
C\|f\|_{L^2(\R^d)}\, .
\end{equation*} 
\end{proposition}
\begin{proof}
Set 
$$
u=e^{-itP} f\,.
$$
Using Proposition~\ref{LP}, we can write
\begin{equation}\label{w0}
\|(1-\chi)u\|_{L^{p}_{T}L^{q}}\leq
C\|f\|_{L^2}+
C\big( \sum_{h^{-1}\,:\,\,{\rm dyadic}}\,
\|\varphi(h^2 P)(1-\chi) u\|^{2}_{L^{p}_{T}L^q}\big)^{\frac{1}{2}}\, .
\end{equation}
Let $\varphi_1\in C_{0}^{\infty}(\R\backslash \{0\})$ be equal to one on the
support on $\varphi$. 
We can write
\begin{eqnarray*}
\varphi(h^2 P)(1-\chi) & = & \varphi(h^2 P)\varphi_1(h^2 P)(1-\chi)
\\
& = &
\varphi(h^2 P)(1-\chi)\varphi_1(h^2 P)
+
\varphi(h^2 P)[\chi,\varphi_1(h^2 P)]\,.
\end{eqnarray*}
Using Proposition~\ref{phi} and Proposition~\ref{est1prelim} (see also Remark~\ref{localized}), we have
\begin{equation}\label{w1}
\|\varphi(h^2 P)(1-\chi)\varphi_1(h^2 P)u\|_{L^{p}_{T}L^{q}}\leq
C\|\varphi_1(h^2 P) f\|_{L^2}\,.
\end{equation}
Using the  Schur lemma, we get
$$
\|[\chi,\varphi_1(h^2 P)]\|_{L^2\rightarrow L^2}\leq Ch,\quad\Big\|\big[[\chi,\varphi_1(h^2 P)],\varphi(h^2 P)\big]\Big\|_{L^2\rightarrow L^2}\leq Ch^2\,.
$$
Thus, using Proposition~\ref{phi}, we can write
\begin{eqnarray*}
\|\varphi(h^2 P)[\chi,\varphi_1(h^2 P)]u\|_{L^p_{T}L^q}
& \leq & 
Ch^{-1}\|\varphi(h^2 P)[\chi,\varphi_1(h^2 P)]u\|_{L^p_{T}L^2}
\\
& \leq & Ch\|u\|_{L^p_{T}L^2}+C\|\varphi(h^2 P)f\|_{L^2}
\\
& \leq & 
Ch\|f\|_{L^2}+C\|\varphi(h^2 P)f\|_{L^2} \,.
\end{eqnarray*}
In view of (\ref{w1}) and the last estimate, coming back to (\ref{w0})
completes the proof of Proposition~\ref{est1prelimbis}
\end{proof}
\section{Semi-classical time estimates and applications}
In this section $g$ is a metric satisfying (\ref{h1}), (\ref{h2}).
The next proposition describes the WKB approximation for solutions of the
semi-classical Schr\"odinger equation for times which are small but independent of the
semi-classical parameter. 
This construction is well known (see e.g. \cite{Robert}). 
Here is the precise statement.
\begin{proposition}\label{WKB}
Let $a(x,\xi)$ be a smooth function on $\R^{2d}$ satisfying
\begin{equation}\label{Hyp1}
\big|\partial_{x}^{\alpha}\partial_{\xi}^{\beta}a(x,\xi)\big|\leq
C_{\alpha\beta},\quad \forall\, (\alpha,\beta)\in \N^{2d}\, ,
\end{equation}
\begin{equation}\label{Hyp2}
\exists\, R>0\, : \, a(x,\xi)=0,\quad {\rm for }\quad |\xi|>R\, .
\end{equation}
Then there exists $\alpha>0$, there exists
$$
S(t,x,\xi)\in C^{\infty}\Big([-\alpha,\alpha]\times\R^d\times \R^d \Big)
$$
and a sequence of smooth functions $a_{j}(t,x,\xi)$, $j\geq 0$ compactly
supported with respect to $\xi$ such that for every $u_0\in L^{2}(\R^d)$, every
$h\in]0,1]$,
every $N\in\N$, the solutions of the problem
$$
(ih\partial_{t}+h^2\Delta_{g})u=0,\quad u|_{t=0}=a(x,hD)u_0
$$
can be represented, for $t\in[-\alpha,\alpha]$ as
$$
u(t,x)=J_{N}(t)u_0+R_{N}(t)u_0,
$$
where $R_{N}(t)$ satisfies
\begin{equation}\label{Rest}
\|R_{N}(t)\|_{L^2(\R^d)\rightarrow H^k(\R^d)}\le Ch^{N+1-k},\quad 0\leq k\leq N+1
\end{equation}
uniformly with respect to $t\in [-\alpha,\alpha]$, and
$$
J_{N}(t)u_0=
(2\pi h)^{-d}
\int_{\R^d}
e^{ih^{-1}S(t,x,\xi)}
\Big[\sum_{j=0}^{N}a_{j}(t,x,\xi)h^j\Big]
\widehat{u_0}\big(\frac{\xi}{h}\big)d\xi\, .
$$
Moreover $S(t,x,\xi)$ is the solution of the Hamilton-Jacobi equation
$$
\partial_{t}S+g(x)(\nabla_{x}S,\nabla_{x}S)=0,\quad
S|_{t=0}=x\cdot\xi
$$
and $a_{j}(t,x,\xi)$ are solutions of the transport equations
$$
\partial_{t}a_0+2g(x)(\nabla_{x}S,\nabla_{x}a_0)+\Delta_{g}S\, a_{0}=0,\quad
a_{0}(0,x,\xi)=a(x,\xi)
$$
for $j=0$, and
$$
\partial_{t}a_j+2g(x)(\nabla_{x}S,\nabla_{x}a_j)+\Delta_{g}S\, a_{j}=i\,\Delta_{g}(a_{j-1}),\quad
a_{j}(0,x,\xi)=0
$$
for $j\geq 1$. 
Finally, there exists $C_{N}>0$ such that if $u_0\in L^1(\R^d)$ then for every $t\in[-\alpha,\alpha]$, every
$h\in ]0,1]$,
$$
\|J_{N}(t) u_0\|_{L^{\infty}(\R^d)}\leq \frac{C_N}{(|t|h)^{d/2}}\|u_{0}\|_{L^1(\R^d)}\,.
$$
\end{proposition}
\begin{proof}
The proof of Proposition~\ref{WKB} is given in \cite{BGT1} when
$a(x,\xi)$ is supported in a coordinate patch of the cotangent bundle of a compact manifold.
The analysis in the case here is slightly more delicate since the $L^2$ bound
of the remainder is not straightforward as in \cite{BGT1}. However, using that
for $|t|\leq \alpha\ll 1$ one has
$$
\big|
\nabla_{x}\nabla_{\xi}S(t,x,\xi)-{\rm Id}
\big|\ll 1,
$$
we can apply the standard result for $L^2$ boundedness of
FIO (see e.g. \cite{Robert}) from which (\ref{Rest}) follows.
\end{proof}
After the time rescaling $t\mapsto ht$,
as in \cite{BGT1} an application of the Keel-Tao theorem \cite{KT} gives
the following Strichartz inequalities (homogeneous and non homogeneous) 
for the Schr\"odinger equation on semi-classical time intervals.
\begin{proposition}\label{semi}
Let $\varphi\in C_{0}^{\infty}(\R)$. 
Then there exist $\alpha>0$ and $C>0$ such that
for every interval $J\subset\R$ of size $\leq\alpha h$, $h\in]0,1]$, if $u$
solves
$$
(i\partial_{t}-P)u=0,\quad u|_{t=0}=\varphi(h^2P)u_0,\quad
u_0\in L^2(\R^d)
$$
then
$$
\|u\|_{L^p(J;L^q(\R^d))}\leq
C\|\varphi(h^2P)u_0\|_{L^2(\R^d)}\, ,
$$
provided $(p,q)$ satisfies (\ref{adm}).
\\
Moreover, if $u$ solves
\begin{equation}\label{nh}
(i\partial_{t}-P)u=\varphi(h^2P)f,\quad 
u|_{t=0}=0
\end{equation}
then
\begin{equation}\label{p1}
\|u\|_{L^p(J;L^q(\R^d))}\leq
C\|\varphi(h^2P)f\|_{L^{p_1}(J;L^{q_1}(\R^d))}\, ,
\end{equation}
provided  $(p,q)$ and $(\frac{p_1}{p_1-1},\frac{q_1}{q_1-1})$ satisfy (\ref{adm}).
\end{proposition}
As in \cite{BGT1}, Proposition~\ref{semi} yields a Strichartz inequality, with
derivative loss in classical Sobolev spaces.
\begin{proposition}\label{perte}
Let $T>0$. Then there exists $C>0$ such that if $u$ solves
$$
(i\partial_{t}-P)u=0,\quad u|_{t=0}=u_0,\quad
u_0\in H^{\frac{1}{p}}(\R^d)
$$
then
$$
\|u\|_{L^p_{T}L^q(\R^d)}\leq
C\|u_0\|_{H^{\frac{1}{p}}(\R^d)}\, ,
$$
provided $(p,q)$ satisfies (\ref{adm}).
\end{proposition}
\begin{proof}
Thanks to Proposition~\ref{LP}, it suffices to prove that for every
$\varphi\in C_{0}^{\infty}(\R)$
there exists $C>0$ such that for every $h\in ]0,1]$, every $f\in L^2(\R^d)$,
$$
\|\exp(-itP)\varphi(h^2P)f\|_{L^p_{T}L^q}\leq
Ch^{-\frac{1}{p}}\|f\|_{L^2}\, .
$$
We split the interval $[-T,T]$ into $Ch^{-1}$ intervals of size $\alpha h$,
where $\alpha$ is the real number involved in the statement of
Proposition~\ref{WKB}.
Using the $L^2$ boundedness of $\exp(-itP)$, and applying, 
according to the above splitting,
about $Ch^{-1}$ times Proposition~\ref{semi} yields
$$
\|\exp(-itP)\varphi(h^2P)f\|_{L^p_{T}L^q}^{p}
\leq Ch^{-1}\|f\|_{L^2}^{p}
$$
which completes the proof of Proposition~\ref{perte}.
\end{proof}
The proof of Theorem~\ref{thm1} is now complete.
\\

Next, we state a non homogeneous extension of Proposition~\ref{semi}.
\begin{proposition}\label{nonhom}
Let $\varphi\in C_{0}^{\infty}(\R)$ and $T>0$. Then there exists $C>0$ such
that if $u$ is a solution of
$$
iu_t-Pu=\varphi(h^2P)f,\quad h\in]0,1],\quad f\in
L^2([0,T]\times\R^d)
$$
with initial data
$$
u|_{t=0}=\varphi(h^2P)u_0,\quad u_0\in L^2(\R^d)
$$
then
\begin{equation}\label{tab}
\|u\|_{L^p_{T}L^q}\leq C\|u\|_{L^{\infty}_{T}L^2}+
Ch^{-1/2}\|u\|_{L^{2}_{T}L^2}+Ch^{1/2}\|\varphi(h^2P)f\|_{L^2_{T}L^2},
\end{equation}
provided $(p,q)$ satisfies (\ref{adm}).
\\
Moreover, if $d=3$, the following estimate holds
\begin{equation}\label{tabbis}
\|u\|_{L^2_{T}L^6}\leq C\|u\|_{L^{\infty}_{T}L^2}+
Ch^{-1/2}\|u\|_{L^{2}_{T}L^2}+C\|\varphi(h^2P)f\|_{L^2_{T}L^{6/5}}\, .
\end{equation}
\end{proposition}
\begin{remarque}
By taking $f=0$ and $T\sim h$, we observe that Proposition~\ref{semi} is a
particular case of (\ref{tab}).
\end{remarque}
\begin{remarque}
For $d\geq 4$, an estimate analogous to (\ref{tabbis}) holds. More precisely,
one has to replace $6$ by $\frac{2d}{d-2}$ and $6/5$ by $\frac{2d}{d+2}$.
\end{remarque}
\begin{proof}[Proof of Proposition~\ref{nonhom}]
Clearly we can restrict our considerations to the interval $[0,T]$,
the analysis on $[-T,0]$ being analogous (the sign of $t$ in this discussion
is harmless).
Write
\begin{equation}\label{Duhamel}
u(t)=e^{-itP}\varphi(h^2P)u_0-i\int_{0}^{t}e^{-i(t-\tau)P}\varphi(h^2P)f(\tau)d\tau\,.
\end{equation}
Observe that $u\in C([-T,T];L^2(\R^d))$. 
If $T\leq\alpha h$, then Proposition~\ref{semi} and the triangle inequality
give
\begin{eqnarray*}
\|u\|_{L^p_{T}L^q} & \leq & C\|\varphi(h^2P)u_0\|_{L^2}+
C\|\varphi(h^2P)f\|_{L^1_{T}L^2}
\\
& \leq &
C\|u\|_{L^{\infty}_{T}L^2}+Ch^{1/2}\|\varphi(h^2P)f\|_{L^2_{T}L^2},
\end{eqnarray*}
provided $(p,q)$ satisfies (\ref{adm}). Hence we can suppose that $T\geq\alpha
h$. Consider a splitting of $[0,T]$ :
$$
[0,T]=[0,a]\cup J_1\cup\dots\cup J_k\cup [b,T],
$$
where for $j=1,\dots k$, there exists $c_j$ such that
$$
J_{j}=\big[c_j-\frac{\alpha h}{8},c_j+\frac{\alpha h}{8}\big]\subset
\big[\frac{\alpha h}{8},T-\frac{\alpha h}{8}\big]\, .
$$
We also suppose that $a\leq \alpha h$ and $T-b\leq\alpha h$. We may also
suppose that
$$
c_1<c_2<\dots<c_k\, .
$$
Clearly $k\lesssim h^{-1}$.
Coming back to the Duhamel formula (\ref{Duhamel}), using
Proposition~\ref{semi}, we obtain that for $(p,q)$ satisfying (\ref{adm}),
\begin{equation}\label{I}
\|u\|_{L^{p}([0,a];L^q)}\leq
C\|u\|_{L^{\infty}_{T}L^2}+
Ch^{1/2}\|\varphi(h^2P)f\|_{L^2([0,a]\times \R^d)}\, .
\end{equation}
Similarly, we estimate the contribution of $[b,T]$ by writing,
\begin{equation}\label{II}
\|u\|_{L^{p}([b,T];L^q)}\leq
C\|u\|_{L^{\infty}_{T}L^2}+
Ch^{1/2}\|\varphi(h^2P)f\|_{L^2([b,T]\times \R^d)}\, .
\end{equation}
We next define the intervals
$$
J_{j}':=J_{j}+\big[-\frac{\alpha h}{8},\frac{\alpha h}{8}\big],
\quad j=1,\dots k\, .
$$
Observe that $J_{j}'\subset [0,T]$. Let us fix $\psi\in C_{0}^{\infty}(\R)$
such that $\psi=1$ on $[-1/8,1/8]$ and ${\rm supp}\,\psi \subset [-1/4,1/4]$.
Set
$$
\psi_{j}(t):=\psi\big(\frac{t-c_j}{\alpha h}\big),\quad j=1,\dots k
$$
and
$$
u_{j}(t):=\psi_{j}(t)u(t)\, .
$$
Observe that $u_{j}(t)=u(t)$ for $t\in J_j$ and that $u_j$ solves the equation
$$
(i\partial_t-P)u_j=i\psi'_{j}\,u+\psi_{j}\,\varphi(h^2P)f,\quad
u_{j}(0)=0\, .
$$
Hence, by writing the Duhamel formula for $u_j$, using
Proposition~\ref{semi} and the triangle inequality, we get for $(p,q)$
satisfying (\ref{adm}),
\begin{eqnarray*}
\|u\|_{L^p(J_j,L^q)} & \leq & \|u_j\|_{L^p(J'_{j};L^q)}
\\
& \leq &
Ch^{-1}\|u\|_{L^1(J'_{j};L^2)}+C\|\varphi(h^2P)f\|_{L^1(J'_{j};L^2)}
\\
& \leq &
Ch^{-1/2}\|u\|_{L^2(J'_{j}\times\R^d)}+Ch^{1/2}\|\varphi(h^2P)f\|_{L^2(J'_{j}\times\R^d)}\, .
\end{eqnarray*}
Summing over $j=1,\dots k$, since $p\geq 2$, we get
\begin{multline*}
\|u\|_{L^p([a,b];L^q)}^{p} \leq Ch^{-\frac{p}{2}}\sum_{j=1}^{k}\|u\|_{L^2(J'_{j}\times \R^d)}^{p}
+Ch^{\frac{p}{2}}\sum_{j=1}^{k}\|\varphi(h^2P)f\|_{L^2(J'_{j}\times \R^d)}^{p}
\\
\leq  Ch^{-\frac{p}{2}}
\Big(\sum_{j=1}^{k}\|u\|_{L^2(J'_{j}\times \R^d)}^{2}\Big)^{\frac{p}{2}}
+Ch^{\frac{p}{2}}\Big(\sum_{j=1}^{k}\|\varphi(h^2P)f\|_{L^2(J'_{j}\times \R^d)}^{2}
\Big)^{\frac{p}{2}}
\\
\leq Ch^{-\frac{p}{2}}\|u\|_{L^2_{T}L^2}^{p}+Ch^{\frac{p}{2}}\|\varphi(h^2P)f\|_{L^2_{T}L^2}^{p}\,.
\end{multline*}
Hence
$$
\|u\|_{L^p([a,b];L^q)} \leq
Ch^{-\frac{1}{2}}\|u\|_{L^2_{T}L^2}+
Ch^{\frac{1}{2}}\|\varphi(h^2P)f\|_{L^2_{T}L^2}\,.
$$
Coming back to (\ref{I}) and (\ref{II}) completes the proof of (\ref{tab}).

Let us now turn to the proof of (\ref{tabbis}).
It follows along the same lines as the proof of  (\ref{tab}).
Indeed, with the above notations, using (\ref{p1}) with $d=3$ and $p=p_1=2$, we obtain that
\begin{eqnarray*}
\|u\|_{L^2(J_j,L^6)} & \leq & \|u_j\|_{L^2(J'_{j};L^6)}
\\
& \leq &
Ch^{-1}\|u\|_{L^1(J'_{j};L^2)}+C\|\varphi(h^2P)f\|_{L^2(J'_{j};L^{6/5})}
\\
& \leq &
Ch^{-1/2}\|u\|_{L^2(J'_{j}\times\R^d)}+
C\|\varphi(h^2P)f\|_{L^2(J'_{j};L^{6/5})}\, .
\end{eqnarray*}
Squaring and summing over $J_{j}$ gives
$$
\|u\|_{L^{2}([a,b];L^6)}\leq Ch^{-1/2}\|u\|_{L^2_{T}L^2}+C\|\varphi(h^2P)f\|_{L^2_{T}L^{6/5}}\, .
$$
Similar estimates holds on $[0,a]$ and $[b,T]$ which ends the proof of (\ref{tabbis}).
This completes the proof of Proposition~\ref{nonhom}.
\end{proof}
\section{Using the non trapping assumption}
The proof of Theorem~\ref{thm2} will be completed, once we prove the following statement.
\begin{theoreme}\label{thm3}
Let the metric $g$ be non trapping and satisfying (\ref{H1}) and (\ref{H1bis}).
Then for every $T>0$ and $\chi\in C_{0}^{\infty}(\R^d)$ there exists a
constant $C$ such that if $u$ solves
\begin{equation}\label{zvezda}
iu_t+\Delta_{g}u=0,\quad u|_{t=0}=u_0\in L^2(\R^d)
\end{equation}
then
$$
\|\chi u\|_{L^{p}_{T}L^q(\R^d)}\leq C\|u_0\|_{L^2(\R^d)}
$$
provided the couple $(p,q)$ satisfies the admissibility condition (\ref{adm}).
\end{theoreme}
\begin{proof}
The non trapping assumption is only needed for the next proposition.
\begin{proposition}\label{doi}
The solution of (\ref{zvezda}) satisfies
$$
\|\chi u\|_{L^2_{T}H^{\frac{1}{2}}}\leq C\|u_0\|_{L^2(\R^d)}\,.
$$
\end{proposition}
We refer to \cite{Doi} for a proof of Proposition~\ref{doi}.
Such estimates can be seen as a consequence of the smooth perturbation theory of
Kato (see \cite[Chapter XIII.7]{RS}).
Let us also recall (see e.g. \cite{BGT2}) that, via a quite general argument
using the Fourier transform in $t$, one can freeze the time and
Proposition~\ref{doi} follows from (the time independent) estimates on the
resolvent of $\Delta_{g}$, namely
\begin{equation}\label{burnol}
\|\chi(P-\lambda\pm i0)^{-1}\chi\|_{L^2\rightarrow L^2}\leq C\langle \lambda
\rangle^{-1/2}\, ,\quad \lambda\gg 1.
\end{equation}
Recall that such resolvent estimates were extensively studied 
in particular in connection with the local energy decay for the wave equation
$(\partial_{t}^{2}-\Delta_{g})u=0$.
For a proof of (\ref{burnol}), we refer for instance to \cite{RobertENS,burqIMRN}. 

Notice that Proposition~\ref{doi} is the only place in the proof of
Theorem~\ref{thm2} where we use the non trapping assumption.

Let us now come back to the proof of Theorem~\ref{thm3}. It will be a
suitable combination of Proposition~\ref{doi} and  Proposition~\ref{nonhom}. 
Let us fix $\varphi\in C_{0}^{\infty}(\R\backslash \{0\})$. Set
$$
v(t):=\varphi(h^2P)\, \chi\, u(t)\, .
$$
Then $v$ solves
\begin{equation*}
(i\partial_{t}-P)v=-\varphi(h^2P)[P,\chi]u,\quad
v|_{t=0}=\varphi(h^2P)\chi u_0\, .
\end{equation*}
We can now apply Proposition~\ref{nonhom} to $v$ which gives
\begin{multline*}
\|v\|_{L^p_{T}L^q}\leq C\|\varphi(h^2P)\chi u\|_{L^{\infty}_{T}L^2}+
Ch^{-\frac{1}{2}}\|\varphi(h^2P)\chi u\|_{L^{2}_{T}L^2}
\\
+Ch^{\frac{1}{2}}\|\varphi(h^2P)[P,\chi] u\|_{L^{2}_{T}L^2}
:=Q_1+Q_2+Q_3,
\end{multline*}
provided $(p,q)$ is satisfying (\ref{adm}). 
We now estimate separately $Q_1$, $Q_2$ and $Q_3$.

{\bf Bound for $Q_1$. } Using the functional calculus and the Schur lemma, we
get
\begin{equation}\label{Sh}
\Big\|\big[\varphi(h^2P),\chi\big] w\Big\|_{L^2(\R^d)}\leq Ch\|w\|_{L^2(\R^d)}\, .
\end{equation}
This implies that
\begin{eqnarray*}
Q_1 & \leq & C\|\chi\, \varphi(h^2P)\,
u\|_{L^{\infty}_{T}L^2}+Ch\|u\|_{L^{\infty}_{T}L^2}
\\
& \leq &
C\| \varphi(h^2P)\,u\|_{L^{\infty}_{T}L^2}
+Ch\|u\|_{L^{\infty}_{T}L^2}
\\
& = &
C\|\varphi(h^2P)\,u_0\|_{L^2}+Ch\|u_0\|_{L^2},
\end{eqnarray*}
where in the last line we used that $\exp(-itP)$ is an $L^2$
isometry.

{\bf Bound for $Q_2$.} Let $\tilde{\varphi}\in
C_{0}^{\infty}(\R\backslash\{0\})$ which is equal to one on the support of
$\varphi$. Then using (\ref{Sh}), we get
$$
Q_2\leq 
Ch^{-\frac{1}{2}}\|\varphi(h^2P)\chi
\tilde{\varphi}(h^2P)u\|_{L^{2}_{T}L^2}
+Ch^{\frac{1}{2}}\|u\|_{L^{\infty}_{T}L^2}:=Q_{21}+Q_{22}\, .
$$
Since the support of $\varphi$ does not meet the origin, we can use
(\ref{koko}) and thus
$$
Q_{21}\leq C\|\chi
\tilde{\varphi}(h^2P)u\|_{L^2_{T}H^{\frac{1}{2}}(\R^d)}\, .
$$
An application of Proposition~\ref{doi} gives
$$
Q_{21}\leq C\|\tilde{\varphi}(h^2P)u_0\|_{L^2}\,\, .
$$
Clearly
$$
Q_{22}=Ch^{\frac{1}{2}}\|u_0\|_{L^2}\, .
$$
Thus
$$
Q_2\leq  C\|\tilde{\varphi}(h^2P)u_0\|_{L^2}+Ch^{\frac{1}{2}}\|u_0\|_{L^2}\, .
$$
{\bf Bound for $Q_3$.}
Let us take again $\tilde{\varphi}\in
C_{0}^{\infty}(\R\backslash\{0\})$ which equals one on the support of
$\varphi$. An application of Schur lemma yields the bound
$$
\Big\|
\big[
\tilde{\varphi}(h^2P),[\chi,P]
\big]w
\Big\|_{L^{2}(\R^d)}\leq C\|w\|_{L^2(\R^d)}\, .
$$
Therefore
$$
Q_3\leq
Ch^{\frac{1}{2}}\|\varphi(h^2P)[\chi,P]\tilde{\varphi}(h^2P)u\|_{L^2_{T}L^2}
+Ch^{\frac{1}{2}}\|u_0\|_{L^2}:=Q_{31}+Q_{32}\, .
$$
Let us next fix a $\tilde{\chi}\in C_{0}^{\infty}(\R^d)$ which is equal to one
on the support of $\chi$. Then
$$
Q_{31}=
Ch^{\frac{1}{2}}\|\varphi(h^2P)
[\chi,P]\tilde{\chi}\tilde{\varphi}(h^2P)u\|_{L^2_{T}L^2}\,.
$$
Using (\ref{koko}) with $s=-1/2$, we obtain the bound
$$
Q_{31}\leq C
\|[\chi,P]\tilde{\chi}\tilde{\varphi}(h^2P)u\|_{L^2_{T}H^{-\frac{1}{2}}}
\leq
C \|\tilde{\chi}\tilde{\varphi}(h^2P)u\|_{L^2_{T}H^{\frac{1}{2}}}
\,,
$$
where in the last line we used that $[\chi,P]$ is a first order
differential operator with $C_{0}^{\infty}(\R^d)$ coefficients.
A use of Proposition~\ref{doi} yields
$$
Q_{31}\leq C\|\tilde{\varphi}(h^2P)u_0\|_{L^2}
$$
and therefore
$$
Q_{3}\leq
C\|\tilde{\varphi}(h^2P)u_0\|_{L^2}+Ch^{\frac{1}{2}}\|u_0\|_{L^2}\, .
$$
Using the above bounds for $Q_1, Q_2, Q_3$, we arrive at the bound
\begin{equation}\label{va}
\|\varphi(h^2P)\chi u\|_{L^{p}_{T}L^q}\leq 
C\|\tilde{\varphi}(h^2P)u_0\|_{L^2}
+
Ch^{\frac{1}{2}}\|u_0\|_{L^2}\, ,
\end{equation}
provided $(p,q)$ is satisfying (\ref{adm}). With (\ref{va}) in hand the proof
of Theorem~\ref{thm3} is reduced to an application of the Littlewood-Paley
square function theorem. Indeed, consider a Littlewood-Paley partition of the
identity
\begin{equation}\label{identity}
{\rm Id}\,=
\varphi_{1}(P)+\sum_{h^{-1}\,:\,\,{\rm dyadic}}\,\varphi(h^2P),
\end{equation}
where $\varphi_{1}\in C_{0}^{\infty}(\R)$, $\varphi\in
C_{0}^{\infty}(\R\backslash\{0\})$.
Using Proposition~\ref{LP}, we obtain that
$$
\|\chi u\|_{L^{p}_{T}L^q}\leq
C\|u_0\|_{L^2}+C\big( \sum_{h^{-1}\,:\,\,{\rm dyadic}}\,
\|\varphi(h^2P)\chi u\|^{2}_{L^{p}_{T}L^q}\big)^{\frac{1}{2}}\, .
$$
Coming to the crucial bound (\ref{va}), using that $\tilde{\varphi}\in
C_{0}^{\infty}(\R\backslash\{0\})$
for the first term in the right hand-side of (\ref{va}), and, summing
geometric series for the second, we arrive at the bound
$$
\|\chi u\| _{L^{p}_{T}L^q}\leq C\|u_0\|_{L^2}\,\, .
$$
This completes the proof of Theorem~\ref{thm3}.
\end{proof}
\begin{remarque}\label{oops}
Let us observe that the constants depending on the time intervals $[0,T]$
in all statements in this paper remain bounded as $T$ varies within a compact
set. In other words the only possible blow up of these
constants is as $T\rightarrow\infty$.
\end{remarque}
\section{Non homogeneous estimates and nonlinear applications}
The aim of this section is to give some applications of the estimates
established in the previous sections to the Nonlinear Schr\"odinger equation
\begin{equation}\label{NLS}
(i\partial_{t}+\Delta_{g})u=F(u),\quad u|_{t=0}=u_{0},
\end{equation}
where $u(t):\R^d\rightarrow \C$. The function $F(z)$, $z\in\C$, is assumed to 
be smooth and vanishing at $z=0$. Moreover, we suppose that
$F=\bar{\partial}V$ with a real valued ``potential'' $V$ satisfying the gauge invariance
assumption
$$
V(\omega z)=V(z),\quad \forall\, \omega\in S^1,\,\, \forall\, z\in\C\, .
$$
In addition, we suppose that for some $\alpha>1$,
\begin{equation}\label{rast}
\big|
\partial^{k_1}\bar{\partial}^{k_2}V(z)
\big|\leq C_{k_1,k_2}\langle z\rangle^{1+\alpha-k_1-k_2}\,.
\end{equation}
The real number $\alpha$ involved in (\ref{rast}) corresponds to the ``degree'' of
the nonlinear interaction. The problem (\ref{NLS}) may, at last formally, be
seen as a Hamiltonian PDE in an infinite dimensional phase space, with
Hamiltonian
\begin{equation}\label{energy}
H(u,\bar{u})=\int_{\R^d}|\nabla_{g}u|^{2}+\int_{\R^d}V(u)
\end{equation}
and canonical coordinates $(u,\bar{u})$
(in (\ref{energy}) we integrate with respect to the volume element associated
to $g$). 
Therefore the quantity (\ref{energy})
is formally conserved by the flow of (\ref{NLS}). Another formally conserved
quantity by the flow of (\ref{NLS}) is the $L^2$ norm of $u$. In this section
we make the defocusing assumption
$$
V(z)\geq 0
$$
on the potential $V$. Under this assumption the $H^1(\R^d)$ norm of the
solutions of (\ref{NLS}) may be expected to be controlled uniformly in time
under the evolution of (\ref{NLS}). Therefore the study of (\ref{NLS}) in the
space $H^1(\R^d)$ is of particular interest.

In the study of (\ref{NLS}), $L^p$ analogues, $1<p<+\infty$, 
of (\ref{lesno}) are useful.
More precisely one has the bounds
\begin{equation}\label{lesno-bis}
C_{s,p}^{-1}\|(P+1)^{s/2}u\|_{L^p(\R^d)}\leq \|u\|_{W^{s,p}(\R^d)}\leq
C_{s,p}\|(P+1)^{s/2}u\|_{L^p(\R^d)}.
\end{equation}
where $1<p<+\infty$.
Estimates (\ref{lesno-bis}) follow from the $L^p$, $1<p<+\infty$ boundedness
of zero order pseudo differential operators (see e.g. \cite{Sogge}). 

In this section, we give the rather standard consequences of
Theorem~\ref{thm1} and Theorem~\ref{thm2} to the $H^1$ theory for (\ref{NLS}).
We start with  a general result in dimension two.
\begin{theoreme}\label{thm4}
Let $\alpha>1$ be an arbitrary real number and let $g$ be a metric on $\R^2$
satisfying (\ref{h1}), (\ref{h2}). Then for every $u_0\in H^1(\R^2)$ there
exists a unique global solution $u\in C(\R;H^1(\R^2))$ of (\ref{NLS}).
\end{theoreme}
\begin{proof}
Using Theorem~\ref{thm1}, we obtain the estimate
\begin{equation}\label{trivial}
\|\exp(it\Delta_{g})u_0\|_{L^p_{T}L^q}\leq C\|u_0\|_{H^{\frac{1}{p}}}\,,
\end{equation}
provided $(p,q)$ is satisfying (\ref{adm}). With (\ref{trivial}) in hand, the
proof of Theorem~\ref{thm4} consists in word by word repetition of the
analysis in \cite[sec. 3.1 and 3.2]{BGT1}.
\end{proof}
Without the non trapping assumption, in dimension three, one can only get the following global
existence result.
\begin{theoreme}\label{thm5}
Consider the cubic defocusing NLS
\begin{equation}\label{cubic}
(i\partial_{t}+\Delta_{g})u=|u|^2 u,\quad u|_{t=0}=u_{0},
\end{equation}
where $g$ is a metric on $\R^3$ satisfying (\ref{h1}), (\ref{h2}).
Let $s>1$. Then for every $u_0\in H^s(\R^3)$ there
exists a unique global solution $u\in C(\R;H^s(\R^3))$ of the Cauchy problem 
(\ref{cubic}).
\end{theoreme}
\begin{proof} It relies on Lemmas \ref{exun} and \ref{boot} below. The first one is a
  {\it local} existence and uniqueness result which  is proved in \cite{BGT1} 
(using the 3 dimensional analogues of (\ref{trivial}) which follow from
Theorem \ref{thm1}). 
\begin{lemme} \label{exun} Let $ s > 1 $ and $ p > 2 $ be 
such that $  s > \frac{3}{2} -
  \frac{1}{p}  $. Set $ \sigma = s - 1/p $ and let $ q > 2 $ be 
such that $ (p,q)  $ satisfy (\ref{adm}). 
Then, for all $ u_0 \in H^s $ and all $ t_0 \in \R
  $, there exists  $ \epsilon > 0 $ and a unique  
$$ u \in C ([t_0- \epsilon , t_0 + \epsilon ],H^s) \cap L^p ( [t_0- \epsilon ,
t_0 + \epsilon ] , W^{\sigma,q} )  $$
such that
\begin{eqnarray}
 u (t) = e^{i(t-t_0)\Delta_g} u_0 - i \int_0^t e^{i(t-\tau)\Delta_g}
|u(\tau)|^2 u (\tau) d \tau,    \label{eitr}
\end{eqnarray}
for all $  t \in [t_0- \epsilon , t_0 +
\epsilon ]  $
\end{lemme}
The key step to get a global existence result is given by the following 
statement.
\begin{lemme} \label{boot}  With $s,p,\sigma,q $ as in Lemma \ref{exun}, 
the following holds true: 
if there exists  $ T > 0 $ and 
$$ u \in \bigcup_{0 \leq T^{\prime} < T} C ([0 , T^{\prime} ],H^s) \cap L^p ( [ 0
, T^{\prime} ] , W^{\sigma,q} ) $$
solution of (\ref{eitr}) (with $t_0 = 0$) 
on $ [0,T^{\prime}] $ for all $T^{\prime}$
 such that $0 \leq  T^{\prime} < T
$, then there exists $ C $ such that
$$ \sup_{[0,T^{\prime}]}||u(t)||_{H^s} + \int_0^{T^{\prime}}
||u(t)||^2_{L^{\infty}} d t \leq C $$
for all $ 0 \leq T^{\prime} < T $.
\end{lemme}
Before proving this lemma, let us show how we obtain Theorem \ref{thm5}. We
consider $ T := \sup \{ T^{\prime} > 0 \ | \ u \  \mbox{solves} \ (\ref{eitr})
\ \mbox{(with} \  t_0 = 0)  \mbox{ on} \ [0,T^{\prime}] \} $ 
and we argue by contradiction, assuming that $ T < \infty $. 
Indeed, using standard non linear estimates and Corollary 2.10 of \cite{BGT1},
Lemma \ref{boot} shows that $ u(T):= \lim_{t \rightarrow T} u(t) $ exists in $ H^s $
and that $ u \in L^p ([0,T],W^{\sigma,q}) $. Thus, by Lemma \ref{exun}, we can continue the
solution on $ [T,T+\epsilon] $, for some $ \epsilon > 0 $, with initial data $
u (T) $ which  yields a
contradiction. 
Of course, we argue similarly for negative times and 
this proves the existence and the uniqueness of a solution in  
$ C ([-T , T ],H^s) \cap L^p ( [ -T , T ] , W^{\sigma,q} ) $ for all $ T > 0 $, 
hence the existence of a solution in $ C (\R ,H^s) $. We omit the proof of the
 uniqueness in $ C (\R,H^s)  $ since it follows as in \cite[3.2]{BGT1} 
and since, here, the main point is the global existence.
Let us finally notice that uniqueness of weak $H^1$ solutions can be
established as in \cite{BGT1}.
\end{proof}

We now turn to the proof of Lemma \ref{boot}.
\begin{proof} The first tool comes from the
conservation laws which imply that there is a constant $C$ independent of $t$
(only depending on $\|u_0\|_{H^1}$)
such that
\begin{equation}\label{H1control}
\|u(t)\|_{H^1(\R^3)}\leq C
\end{equation}
as far as the solution exists, {\it i.e.} on $ [ 0 ,T ) $ here.
The rigorous justification of these conservation laws requires a standard approximation argument (see e.g. \cite{GiDEA})
The  key quantity in this discussion is $\|u\|_{L^{2}_{T}L^{\infty}}$ (the
number $2$ is reflecting the cubic nature of the nonlinearity).
Consider again the Littlewood-Paley partition of the identity
(\ref{identity}).
Then $v:=\varphi(h^2P)u$ solves the problem
$$
(i\partial_{t}-P)v=\varphi(h^2P)\big(|u|^2 u\big)\, .
$$
Using Proposition~\ref{phi}, Proposition~\ref{nonhom} and the bound
(\ref{H1control}), we obtain that for all $ \theta \leq \inf (1,T/2)  $ (see Remark~\ref{oops})
\begin{multline*}
\|\varphi(h^2P)u\|_{L^{2}_{\theta}L^{\infty}}
\leq  Ch^{-1/2}
\|\varphi(h^2P)u\|_{L^{2}_{\theta}L^{6}}
\\
\leq 
C h^{1/2}||u||_{L^{\infty}_{\theta} H^1} +Ch^{-1}\|\varphi(h^2P)u\|_{L^{2}_{\theta}L^{2}}
+Ch^{-1/2}\|\varphi(h^2P)(|u|^2 u)\|_{L^2_{\theta}L^{6/5}}\,.
\end{multline*}
Next, using Proposition~\ref{inverse-bis}, the Sobolev inequality and (\ref{H1control}), we obtain that
\begin{multline*}
h^{-1/2}\|\varphi(h^2P)(|u|^2 u)\|_{L^2_{\theta}L^{6/5}}\leq
Ch^{1/2}\||u|^2 u\|_{L^2_{\theta}W^{1,6/5}}
\leq
\\
\leq Ch^{1/2}\|u\|_{L^2_{\theta}H^1}\leq C\sqrt{h \theta}\,.
\end{multline*}
In summary,
$$
\|\varphi(h^2P)u\|_{L^{2}_{\theta}L^{\infty}}\leq
C h^{1/2} +C\sqrt{h \theta}+Ch^{-1}\|\varphi(h^2P)u\|_{L^{2}_{\theta}L^{2}}\,.
$$
Next, for $N\in\N$, using the Cauchy-Schwarz inequality we get
$$
\sum_{h^{-1}\leq N}\, h^{-1}\|\varphi(h^2P)u\|_{L^{2}_{\theta}L^{2}}\leq
C\sqrt{\log N}\|u\|_{L^2_{\theta}H^1}\leq C\sqrt{\theta\log N}\,.
$$
On the other hand, since $s>1$, we estimate the high frequencies as follows
$$
\sum_{h^{-1}\geq N}\, h^{-1}\|\varphi(h^2P)u\|_{L^{2}_{\theta}L^{2}}\leq
CN^{-(s-1)}\|u\|_{L^2_{\theta}H^s}\,.
$$
By taking $N\approx (2+\|u\|_{L^{\infty}_{\theta}H^s})^{\frac{1}{s-1}}$, we deduce that
\begin{eqnarray} \label{ldeux}
\|u\|_{L^2_{\theta}L^{\infty}}\leq C+C\big[
\theta(1+\log(2+\|u\|_{L^{\infty}_{\theta}H^s}))
\big]^{\frac{1}{2}}\,.
\end{eqnarray}
Coming back to the integral equation (\ref{eitr}) and using the Gronwall lemma, we
obtain that
\begin{equation}\label{key}
\|u\|_{L^{\infty}_{\theta}H^s}\leq \|u_0\|_{H^s}\,\,e^{C\|u\|^{2}_{L^2_{\theta}L^{\infty}}}
\leq C\|u_0\|_{H^s}\big[
2+\|u\|_{L^{\infty}_{\theta}H^s}
\big]^{\Lambda \theta},
\end{equation}
where $\Lambda$ is a real number depending only on the a priori bound (\ref{H1control}).
Therefore, if we take $\theta$ such that $\Lambda \theta \leq 1/2$, 
we obtain that
\begin{eqnarray}
||u||_{L^{\infty}_{\theta}H^s} \leq C \left( ||u_0 ||_{H^s} + ||u_0 ||_{H^s}^2
\right) \label{bgt327} .
\end{eqnarray}
 Iterating finitely many times ($ \approx  T / \theta  $ times)
(\ref{bgt327}) and (\ref{ldeux})  yields the result.
\end{proof} 
\begin{remarque}
Once we know that we have a global solution, we  can
control the growth of $ ||u(t) ||_{H^s} $ as $ t \rightarrow  \infty $ since,
by iterations of (\ref{bgt327}), one can easily check that
\begin{equation}
\|u(t)\|_{H^s}\leq C \exp( C \exp( C |t|)) , \qquad t \in \R .  \nonumber
\end{equation}
\end{remarque}
Notice that the results of Theorem~\ref{thm4} and Theorem~\ref{thm5} hold
without the long range assumption (\ref{H1bis}). Moreover, we do not suppose
that the metric is non trapping. If we assume these two conditions, we can
improve Theorem~\ref{thm5} to nonlinearities of higher degree and even get
global existence results in dimensions $d\geq 4$. For that purpose,  we need the
following non homogeneous Strichartz estimate.
\begin{theoreme}\label{thm6}
Suppose that $g$ is a non trapping metric on $\R^d$ satisfying (\ref{H1}),
(\ref{H1bis}).
Then for every $T>0$ there exists $C>0$ such that that if $u$ solves 
\begin{equation}\label{posledno}
(i\partial_{t}+\Delta_{g})u=f,\quad u|_{t=0}=0
\end{equation}
then
\begin{equation}\label{poslednobis}
\|u\|_{L^p_{T}L^{q}(\R^d)}\leq C\|f\|_{L^{p_1}_{T}L^{q_1}(\R^d)},
\end{equation}
provided $(p,q)$, $(\frac{p_1}{p_1-1},\frac{q_1}{q_1-1})$ are satisfying
(\ref{adm}) and $p\neq 2$, $p_1\neq 2$.
\end{theoreme}
\begin{proof}
The proof is a consequence of the following Christ-Kiselev lemma \cite{CK}.
\begin{lemme}\label{Christ}
Let $T>0$ be a real number. Let $B_1$ and $B_2$ be two Banach spaces. Let
$K(t,s)$ be a locally integrable kernel with values in the bounded operators
from $B_1$ to $B_2$. Let $1\leq p< q\leq\infty$. For $t\in [0,T]$, we set
$$
Af(t)=\int_{0}^{T}K(t,s)f(s)ds\,.
$$
Assume that
$$
\|Af\|_{L^q([0,T];B_1)}\leq C\|f\|_{L^p([0,T];B_2)}\, .
$$
Define the operator $\tilde{A}$ as
$$
\tilde{A}f(t)=\int_{0}^{t}K(t,s)f(s)ds\,,\quad t\in[0,T]\,.
$$
Then there exists $\tilde{C}>0$ such that
$$
\|\tilde{A}f\|_{L^q([0,T];B_1)}\leq \tilde{C}\|f\|_{L^p([0,T];B_2)}\, .
$$
\end{lemme}
We refer to \cite{SS} for a proof of Lemma~\ref{Christ}, in the form stated here.
Let us now return to the proof of Theorem~\ref{thm6}.
The solution of (\ref{posledno}) is given by
$$
u(t)=\int_{0}^{t}\exp(i(t-\tau)\Delta_{g})f(\tau)d\tau\,.
$$
Consider $v(t)$ defined by
$$
v(t)=\int_{0}^{T}\exp(i(t-\tau)\Delta_{g})f(\tau)d\tau\,.
$$
Then
$$
v(t)=\exp(it\Delta_{g})\int_{0}^{T}\exp(-i\tau\Delta_{g})f(\tau)d\tau
$$
and by invoking (\ref{est3}), we get the bound
$$
\|v\|_{L^p_{T}L^q}\leq C\Big\|
\int_{0}^{T}\exp(-i\tau\Delta_{g})f(\tau)d\tau
\Big\|_{L^2}\,.
$$
The dual of (\ref{est3}) yields
\begin{equation}\label{T}
\Big\|\int_{0}^{T}\exp(-i\tau\Delta_{g})f(\tau)d\tau
\Big\|_{L^2}\leq C\|f\|_{L^{\frac{p}{p-1}}_{T}L^{\frac{q}{q-1}}}\, .
\end{equation}
Therefore
$$
\|v\|_{L^p_{T}L^q}\leq C\|f\|_{L^{\frac{p}{p-1}}_{T}L^{\frac{q}{q-1}}}\, .
$$
An application of Lemma~\ref{Christ} gives
$$
\|v\|_{L^p_{T}L^q}\leq C\|f\|_{L^{\frac{p}{p-1}}_{T}L^{\frac{q}{q-1}}}\, ,\quad
p>2.
$$
This proves (\ref{poslednobis}) in the case
$\frac{1}{p}+\frac{1}{p_1}=\frac{1}{q}+\frac{1}{q_1}=1$.
Let us now prove  (\ref{poslednobis}) when $(p,q)$ and $(p_1,q_1)$ are
decoupled.
Since $\exp(it\Delta_g)$ is an isometry on $L^2$, using (\ref{T}) and
Remark~\ref{oops}, we infer that (\ref{poslednobis}) is valid for $(p,q)=(\infty,2)$
and $(p_1,q_1)$ an arbitrary pair satisfying (\ref{adm}). Using the
homogeneous estimate (\ref{est3}) and the Minkowski inequality, we obtain that
(\ref{poslednobis}) is valid for $(p,q)$ an arbitrary pair satisfying
(\ref{adm}) and $(p_1,q_1)=(\infty,2)$. 
Let us now observe that all other cases for $(p,q)$ and $(p_1,q_1)$ in
(\ref{poslednobis}) 
follow from the considered three particular cases by interpolation.
This completes the proof of Theorem~\ref{thm6}.
\end{proof}
It is now a standard and straightforward consequence of Theorem~\ref{thm6},
Theorem~\ref{thm2} and (\ref{lesno-bis})
(see \cite{GVclassic,Kato,Caz,GiDEA}) that one has the following global well-posedness
result for (\ref{NLS}). 
\begin{theoreme}\label{thm7}
Let $d\geq 3$.
Suppose that $\alpha\in ]1,1+\frac{4}{d-2}[$. Let $g$ be a non trapping metric on $\R^d$ satisfying (\ref{H1}),
(\ref{H1bis}). Then for every $u_{0}\in H^1(\R^d)$
there
exists a unique global solution $u\in C(\R;H^1(\R^d))$ of (\ref{NLS}).
\end{theoreme}
\begin{remarque}
Recall that the endpoint Strichartz estimates are not needed for the standard
$H^1$ theory of (\ref{NLS}).
\end{remarque}
\begin{remarque}
Let us emphasize the importance of the non homogeneous Strichartz estimates in
the proof of Theorem~\ref{thm7}. The lack of such estimates under the very
weak hypotheses (\ref{h1}), (\ref{h2}) makes the study of (\ref{NLS}) in
$H^1$ in this case (or in the case of a compact manifold) more difficult and
so far restricted only to small dimensions.
\end{remarque}
\def\cprime{$'$}

\end{document}